\documentclass[11pt]{article}

\usepackage{amsfonts}
\usepackage{graphics}
\usepackage{indentfirst}
\usepackage{color}
\usepackage{cite}
\usepackage{float}
\usepackage{caption}
\usepackage{latexsym}
\usepackage[paper=a4paper, left=1.6cm, right=1.6cm, top=1.8cm, bottom=1.6cm, headheight=5.5pt, footskip=0.8cm, footnotesep=0.8cm, centering, includefoot]{geometry}
\usepackage{amsmath}
\numberwithin{equation}{section}
\allowdisplaybreaks
\usepackage{amssymb}
\usepackage[colorlinks, linkcolor=red]{hyperref}
\hypersetup{urlcolor=red, citecolor=red}
\usepackage[dvips]{epsfig}
\usepackage{amscd}
\usepackage{amsthm}
\usepackage{mathrsfs}
\usepackage{verbatim}
\usepackage[title]{appendix}

\DeclareMathOperator{\loc}{loc}

\newtheorem{theorem}{Theorem}[section]
\newtheorem{remark}{Remark}[section]
\newtheorem{definition}{Definition}[section]
\newtheorem{lemma}{Lemma}[section]

\newtheorem{proposition}{Proposition}[section]

\catcode`@=11 \@addtoreset{equation}{section} \catcode`@=12

\newcommand{\R}{\mathbb{R}}
\newcommand{\RR}{\mathcal{R}}

\title{Scaling-invariant conditions for two-dimensional  compressible Navier--Stokes equations
\thanks{This research was partially supported by National Natural Science Foundation of China (No. 12371227) and Fundamental Research Funds for the Central Universities (No. SWU--KU24001).}
}

\author{Lin Xu,\ Xin Zhong{\thanks{Corresponding author. E-mail addresses: mathxu@email.swu.edu.cn (L. Xu), xzhong1014@amss.ac.cn (X. Zhong).}}
\date{}\\
\footnotesize School of Mathematics and Statistics, Southwest University, Chongqing 400715, P. R. China}
\begin{document}
\maketitle
	
\begin{abstract}
We study the Cauchy problem of two-dimensional isentropic compressible Navier--Stokes equations with either far-field vacuum or non-vacuum. Under suitable smallness conditions on scaling-invariant initial quantities, we establish the global existence and uniqueness of strong solutions involving initial vacuum states and large oscillations. The two far-field regimes present distinct analytical difficulties and require different approaches to obtain the desired {\it a priori} estimates.
		
$\bullet$ For the far-field vacuum case $\rho_\infty=0$, we introduce the scaling-invariant initial quantity
\begin{align*}
\left(\frac{1}{m_0\sqrt{\frac{m_0}{\bar\rho}}}\int_{\mathbb R^2}\rho_0
\sqrt{\frac{m_0}{\bar\rho}+|x|^2}\,dx\right)^{\frac13}
\big[1+(Am_0\bar\rho^\gamma)^{-\frac23}\big]
\bar\rho\big[E_0+m_0\big(\|\nabla u_0\|_{L^2}^2+\|P(\rho_0)\|_{L^2}^2\big)\big],
\end{align*}
where $m_0$ denotes the initial total mass, $\bar\rho$ is a prescribed upper bound for the initial density, and $E_0$ is the initial total energy. The analysis combines mass localization based on the finite first spatial moment with scaling-invariant weighted inequalities and refined time-weighted estimates for the material derivative. These estimates yield short-time and large-time bounds for the effective viscous flux, which are essential to obtain a uniform upper bound of the density.
		
$\bullet$ For the far-field non-vacuum case $\rho_\infty>0$, the corresponding scaling-invariant initial quantity is
\begin{align*}
\rho_\infty^{1-\gamma}\big(\|\nabla u_0\|_{L^2}^2+P_\infty E_0\big),
\end{align*}
where $P_\infty=A\rho_\infty^\gamma$. The proof relies on relative-energy estimates, BMO-based interpolation inequalities, and time-weighted estimates adapted to the positive far-field density.
		
Although the two smallness criteria have different forms, they arise from the same intrinsic scaling of the equations. Our results provide a unified scaling-invariant framework for the global well-posedness of strong solutions in both far-field regimes.
\end{abstract}
	


\setcounter{tocdepth}{2}
\tableofcontents
	
\section{Introduction and main results}
	
The motion of a two-dimensional viscous isentropic compressible fluid in $\mathbb R^2$ is modeled by the Navier--Stokes system
\begin{equation}\label{1}
\begin{cases}
\rho_t+\operatorname{div}(\rho u)=0,\\
\rho(u_t+u\cdot\nabla u)-\mu\Delta u
-(\mu+\lambda)\nabla\operatorname{div}u+\nabla P(\rho)=0.
\end{cases}
\end{equation}
Here $\rho=\rho(x,t)\geq0$ and $u=u(x,t)\in\R^2$ denote the density and velocity, respectively. The pressure and viscosity coefficients satisfy
\begin{align*}
P(\rho)=A\rho^\gamma,\quad A>0,\quad \gamma>1,\quad \mu>0,\quad \mu+\lambda\geq0.
\end{align*}
We consider the Cauchy problem supplemented with the initial data
\begin{align}\label{2}
(\rho,\rho u)(0,x)=\bigl(\rho_{0},\rho_{0}u_{0}\bigr)(x),
\end{align}
and the far-field behavior
\begin{align}\label{3}
\lim_{|x|\rightarrow\infty}(\rho,u)(t,x)=(\rho_{\infty},0),
\end{align}
where $\rho_\infty\ge0$ is a prescribed constant. We treat far-field vacuum case $\rho_\infty=0$ and non-vacuum case $\rho_\infty>0$. In both cases, the initial density may vanish in the interior of the fluid region. Our goal is to establish global strong solutions under assumptions that are exactly invariant under the scaling of \eqref{1}.

\subsection{Background and related results}

Developing a global theory of smooth solutions to the multidimensional compressible Navier--Stokes equations for general large initial data remains an open problem in the mathematical theory of viscous compressible fluids. Classical local solvability for positive density goes back to Nash \cite{Nash1962}, while Matsumura and Nishida \cite{MN1980,MN1983} established global results for sufficiently small perturbations of a equilibrium. Global solutions with lower regular data were developed by Hoff \cite{Hoff1995}. For general large data, the finite-energy weak-solution theory of Lions \cite{Lions1998} and Feireisl, Novotn\'y, and Petzeltov\'a \cite{FNP2001}, together with the developments described in \cite{Feireisl2004}, provides a different framework. Yet uniqueness and global regularity for such weak solutions remain a challenging problem.

The presence of vacuum introduces some additional difficulties since the coefficient of the material derivative of velocity in the momentum equation may vanish. Local strong solvability with vacuum was established by Cho, Choe, and Kim \cite{ChoChoeKim2004}. A major step in the global theory was made by Huang, Li, and Xin \cite{HLX2012}, who obtained three-dimensional global classical solutions with small initial energy involving vacuum and large oscillations, for both far-field vacuum and non-vacuum states. For the two-dimensional Cauchy problem, Luo \cite{Luo2012,Luo2014} developed local and global classical solvability results in the far-field non-vacuum setting, Li and Liang \cite{LiLiang2014} demonstrated a local classical solution permitting vacuum at infinity, and Li and Xin \cite{LiXin2019} proved the global well-posedness and large-time asymptotics of classical solutions with far-field vacuum under a small-energy assumption. These works identify features specific to two-dimensional problem. Spatial weights and the distribution of mass therefore become important when the far-field density vanishes. Meanwhile, large bulk viscosity provides another route to global solutions. Danchin and Mucha \cite{DanchinMucha2017} constructed solutions with large initial velocity in the whole space, exploiting the incompressible dynamics and strong damping of the compressible component. Their subsequent work \cite{DanchinMucha2023} addressed global solutions on the two-dimensional torus with bounded nonnegative under sufficiently large bulk viscosity and suitable control of the initial divergence.

\subsection{Scaling invariance and criticality}

Scaling invariance has played a fundamental role in the development of well-posedness theory for Navier--Stokes equations. The underlying idea is to identify function spaces whose norms remain unchanged under a natural scaling of the system. For the incompressible Navier--Stokes equations
\begin{align*}
\partial_t v-\Delta v+v\cdot\nabla v+\nabla\pi=0,\quad \operatorname{div} v=0,
\end{align*}
the scaling transformation is given by
\begin{align*}
v^\xi(x,t)=\xi v(\xi x,\xi^2t),\quad
\pi^\xi(x,t)=\xi^2\pi(\xi x,\xi^2t),\quad \xi>0.
\end{align*}
In two dimensions, the $L^2$ norm of the initial velocity is invariant under this transformation. More generally, the homogeneous Besov norm $\|v_0\|_{\dot B^s_{p,q}(\R^2)}$ has scaling exponent $1+s-2/p$ so that the critical regularity is $s=2/p-1$. This observation provides a natural connection between the energy structure of the two-dimensional incompressible equations and the critical-space approach initiated by Fujita and Kato \cite{FK64}.

For the isentropic compressible Navier--Stokes equations near a constant equilibrium, the classical critical-space framework is based on the parabolic transformation
\begin{align}\label{g}
\rho^\xi(x,t)=\rho(\xi x,\xi^2t),\quad u^\xi(x,t)=\xi u(\xi x,\xi^2t).
\end{align}
Although the above transformation does not preserve the full system with a fixed pressure law, it determines the critical regularities associated with the transport and viscous terms. In two dimensions, the corresponding critical indices for the density perturbation and velocity are represented by $\dot B^1_{2,1}(\R^2)$ and $\dot B^0_{2,1}(\R^2)$, respectively, with additional low-frequency control required in the global theory. In 2000, Danchin \cite{Danchin2000} established the global well-posedness theory for sufficiently small perturbations of a equilibrium in critical Besov spaces. Subsequently, Danchin \cite{Danchin2005} addressed the uniqueness problem at the critical regularity endpoint in two dimensions. The framework was further extended by Charve and Danchin \cite{CharveDanchin2010} to critical $L^p$-based Besov spaces, allowing highly oscillating initial velocity fields under appropriate smallness assumptions. Further progress in the two-dimensional setting has been made for classes of large initial data. Danchin and Mucha \cite{DanchinMucha2017} established global solutions with arbitrarily large initial velocity and density close to a positive constant, provided that the bulk viscosity is sufficiently large. Later on, Chen and Zhai \cite{ChenZhai2019} generalized their approach to a critical $L^p$ framework involving large highly oscillating velocity fields. More recently, Fujii \cite{Fujii2024} obtained global solutions for large initial data in critical Besov spaces under a sufficiently small Mach number condition. These results demonstrate the effectiveness of critical regularity in the analysis of compressible flows.

However, the perturbative critical-space theories described above generally rely on a strictly positive reference density and suitable smallness or nondegeneracy conditions on the density perturbation. Moreover, their notion of criticality is associated with the parabolic scaling \eqref{g}, under which the pressure law must be rescaled. It should be pointed out that although \eqref{g} preserves the scaling of the transport and viscous terms, it does not preserve the system with a fixed pressure law: its pressure coefficient must be replaced by $\xi^2 A$. Moreover, the standard perturbative critical Besov theory near a equilibrium depends heavily on positivity of the density.

For the polytropic law, there is an exact scaling that keeps $A$, $\gamma$, $\mu$, and $\lambda$ unchanged. Following Lei and Xin \cite{LeiXin2019}, define
\begin{align}\label{2sca}
\rho^\ell(x,t)=\ell^{\frac{2}{\gamma+1}}
\rho \big(\ell x,\ell^{\frac{2\gamma}{\gamma+1}}t\big),\quad u^\ell(x,t)
=\ell^{\frac{\gamma-1}{\gamma+1}}u \big(\ell x,\ell^{\frac{2\gamma}{\gamma+1}}t\big),\quad \ell>0.
\end{align}
Then $(\rho^\ell,u^\ell)$ solves the same system. Recent developments in the global well-posedness theory of compressible Navier--Stokes equations have demonstrated the effectiveness of scaling-invariant smallness conditions. For the non-isentropic system, Li \cite{JinkaiLi2020} established global strong solutions allowing vacuum under a smallness condition invariant under the intrinsic scaling of the equations. For the three-dimensional isentropic compressible Navier--Stokes equations with far-field vacuum, Wen \cite{Wen2025} obtained a global well-posedness result under a scaling-invariant smallness condition. Subsequent developments include the coefficient-coupled smallness condition introduced by Xu and the second author \cite{XZ20261}, the investigation of the role of pressure in critical scaling and low-regularity solutions by the present authors \cite{XuZhong2026}, and our recent work \cite{XuZhong202611} regarding the global well-posedness of strong solutions to the 3D Cauchy problem with a positive far-field density under a scaling-invariant smallness condition. These advances provide further motivation for investigating scaling-invariant global well-posedness in two dimensions. This motivates the current study of both far-field regimes within a unified scaling-invariant framework.

\subsection{Observation and  motivation}
The scaling-invariant approach to the global well-posedness of compressible Navier--Stokes equations has led to important developments in three dimensions. It is therefore natural to ask whether a similar approach can be
developed for the two-dimensional Cauchy problem. However, the transition from three to two dimensions is not a straightforward extension of the existing arguments. One of the main obstacles comes from the different Sobolev structures in two and three dimensions. In three dimensions, the homogeneous Sobolev embedding $\dot H^1(\mathbb R^3)\hookrightarrow L^6(\mathbb R^3)$ provides direct $L^6$ control of the velocity from its gradient. In two dimensions, however, $\|\nabla u\|_{L^2(\mathbb R^2)}$ alone does not control the $L^p$-norm of $u$ for any finite $p$. Moreover, the basic energy inequality only provides control of $\sqrt\rho u$, which cannot be directly converted into an unweighted velocity estimate when vacuum is allowed. Consequently, additional low-frequency information is needed to control the velocity and momentum, as well as the effective viscous flux arising in the density estimates. \textit{An important observation is that the required low-frequency information can be recovered in different ways, depending on the far-field behavior of the density.}

When $\rho_\infty>0$, the positive constant state at infinity provides a natural reference density. The relative potential energy controls the density perturbation, while the weighted Poincar\'e inequality recovers unweighted $L^2$ estimates from density-weighted quantities. Furthermore, the scale-invariant ratio
$\frac{\bar\rho}{\rho_\infty}$ allows the pressure to be controlled in terms of explicit powers of $\rho_\infty$. Indeed, under the density bootstrap assumption $0\le\rho\le2\bar\rho$,
\begin{align*}
\|P(\rho)\|_{L^\infty}\le A(2R)^\gamma\rho_\infty^\gamma,
\end{align*}
and the relative energy provides the corresponding $L^2$ control of $P(\rho)-P_\infty$; see \eqref{eqeq2}. This structure is particularly useful for estimating the nonlinear terms without introducing additional assumptions on the time decay of the pressure.

The situation is fundamentally different when $\rho_\infty=0$. In this case, the relative-energy argument around a positive constant state is no longer available. Instead, we exploit the finite total mass
\begin{align*}
0<m_0\triangleq\int_{\mathbb R^2}\rho_0\,dx<\infty
\end{align*}
and the first spatial moment of the initial density. These quantities provide quantitative information on the spatial distribution of the fluid mass; see \eqref{70} and \eqref{71}. Combined with a weighted Sobolev inequality, this localization property leads to density-weighted $L^p$ estimates, as established in Lemma \ref{zuilem}.

A further observation concerns the intrinsic time scales associated with the two far-field settings. Under the scaling \eqref{2sca}, the total mass in the vacuum case and the far-field density in the non-vacuum case transform according to
\begin{align*}
m_0^\ell=\ell^{-\frac{2\gamma}{\gamma+1}}m_0,\quad
\rho_\infty^\ell =\ell^{\frac{2}{\gamma+1}}\rho_\infty.
\end{align*}
Therefore, both $m_0$ and $\rho_\infty^{-\gamma}$ transform with the same scaling exponent as the time variable.
This suggests two natural scale-adapted time parameters $\frac{t}{m_0}$ and $\rho_\infty^\gamma t$, which lead to the time weights $m_0+t$ and $t(m_0+t)$ in the vacuum case, and $\min\{1,\rho_\infty^\gamma t\}$ in the non-vacuum case; see \eqref{60}--\eqref{61} and \eqref{w0}, respectively.

These observations motivate the introduction of distinct scaling-invariant smallness conditions on the initial data for the two far-field regimes. We develop a unified strategy for establishing global solutions involving the different low-frequency structures associated with the two far-field states.

Before stating our main results, we would like to give some notations that will be used throughout this paper.

\medskip
\noindent
\textbf{Notation:}

$\bullet$ For an integer $k\geq0$ and $1\leq q\leq\infty$, we define
\begin{gather*}
D^{k,q}(\mathbb{R}^{2}) = \{ u \in L^1_{\loc}(\mathbb{R}^{2}) \mid \|\nabla^k u\|_{L^q} < \infty \}, \ \|u\|_{D^{k,q}} = \|\nabla^k u\|_{L^{q}(\mathbb{R}^{2})},\ \|  u \|_{L^{q}}=\|u \|_{L^{q}(\mathbb{R}^{2})}, \\
W^{k,q}(\mathbb{R}^{2}) = L^q(\mathbb{R}^{2}) \cap D^{k,q}(\mathbb{R}^{2}), \ D^k(\mathbb{R}^{2}) = D^{k,2}(\mathbb{R}^{2}), \ H^k(\mathbb{R}^{2}) = W^{k,2}(\mathbb{R}^{2}).
\end{gather*}

$\bullet$
For two nonnegative quantities $a$ and $b$, we write $a\sim b$ if $cb\le a\le Cb$ for some positive constants $c$ and $C$. Unless otherwise specified, $c$ and $C$ denote generic positive constants that may vary from line to line.

$\bullet$
$
{\displaystyle\int} \cdot \, d x \triangleq {\displaystyle\int}_{\mathbb{R}^{2}} \cdot \, d x,\quad  \dot{u}\triangleq u_t+u\cdot\nabla u.
$

$\bullet$
For $R > 0$, $
B_{R} \triangleq\{x \in \mathbb{R}^{2} | | x |<R\}.
$

$\bullet$
$
P_\infty\triangleq A\rho_\infty^\gamma,\quad
p\triangleq P(\rho)-P_\infty,\quad
F\triangleq(2\mu+\lambda)\operatorname{div} u-p,\quad
\omega\triangleq\partial_1u^2-\partial_2u^1.
$

$\bullet$
$\bar\rho\triangleq\|\rho_0\|_{L^\infty},\quad m_0\triangleq{\displaystyle\int}\rho_0\,dx>0,\quad \mathcal{R}\triangleq \frac{\bar\rho}{\rho_\infty}\geq1$.

$\bullet$
$E(t)={\displaystyle\int} \left(\frac{1}{2} \rho \left|u \right|^{2}+G(\rho \mid\rho_\infty)\right)(x,t)\, dx$ is the energy, $E_{0}\triangleq E(0)$ is the initial energy, where $G$ denotes relative potential energy density given by
\begin{align*}
G(\rho\mid\rho_\infty)\triangleq\frac{A}{\gamma-1}\left[\rho^\gamma-\gamma\rho_\infty^{\gamma-1}\rho
+(\gamma-1)\rho_\infty^\gamma\right].
\end{align*}

$\bullet$
$
D(t)\triangleq\|\nabla u(t)\|_{L^{2}}^{2}+\|p(t)\|_{L^{2}}^{2},\quad D_{0}\triangleq D(0),\quad p_{0}\triangleq p(0).
$

We now give the definition of strong solutions considered in this paper.
\begin{definition}[Strong solutions]
Let $T>0$, $q\in(2,\infty)$, $a>1$, and $\eta_0\in(0,1]$. Define
\begin{align*}
\bar{x}\triangleq(e+|x|^2)^{1/2}\log^{1+\eta_0}(e+|x|^2).
\end{align*}
A pair $(\rho\geq0,u)$ is called a strong solution to the Cauchy problem \eqref{1}--\eqref{3} on $\mathbb R^2\times(0,T]$ if it satisfies the equations almost everywhere in $\mathbb R^2\times(0,T)$, attains the prescribed initial data and satisfies the following regularity properties.
	
\medskip
\noindent
{\rm (i)} Far-field vacuum
($\rho_\infty=0$):
\begin{align}\label{def2}
\begin{cases}
\rho\in C([0,T];L^1\cap H^1\cap W^{1,q}),\quad
\bar{x}^{a}\rho\in L^\infty(0,T;L^1\cap H^1\cap W^{1,q}),\\
(\sqrt{\rho}u,\ \nabla u,\ \bar{x}^{-1}u,\ \sqrt{t}\sqrt{\rho}u_t)
\in L^\infty(0,T;L^2),\quad \nabla u\in L^2(0,T;H^1)\cap L^{(q+1)/q}(0,T;W^{1,q}), \\
\sqrt{t}\nabla u\in L^2(0,T;W^{1,q}),\quad (\sqrt{\rho}u_t,\ \sqrt{t}\nabla u_t,\ \sqrt{t}\bar{x}^{-1}u_t)\in L^2(0,T;L^2).
\end{cases}
\end{align}

\medskip
\noindent
{\rm (ii)} Far-field non-vacuum
($\rho_\infty>0$):
\begin{align}\label{def1}
\begin{cases}
(\rho-\rho_\infty)\in C([0,T];H^1\cap W^{1,q}),\quad\rho_t\in L^\infty(0,T;L^2)\cap L^{2q/(q-2)}(0,T;L^q),\\
u\in L^\infty(0,T;H^1)\cap L^2(0,T;H^2)\cap L^{(q+1)/q}(0,T;W^{2,q}), \\
\sqrt{t}\nabla^2u\in L^\infty(0,T;L^2)\cap L^2(0,T;L^q),
\quad\sqrt{t}\nabla u_t\in L^2(0,T;L^2),\\
(\sqrt{\rho}u,\ \sqrt{t}\sqrt{\rho}u_t)\in L^\infty(0,T;L^2),\quad \sqrt{\rho}u_t\in L^2(0,T;L^2).
\end{cases}
\end{align}
\end{definition}

\subsection{Main results}

We consider the Cauchy problem for the problem \eqref{1}--\eqref{3} with $\rho_\infty>0$. In this case, we introduce the following scaling-invariant initial quantity
\begin{align}\label{SS1}
\mathcal S_0^{\mathrm{nvac}}\triangleq\rho_\infty^{1-\gamma}
\left(\|\nabla u_0\|_{L^2}^2+P_\infty E_0\right).
\end{align}
The factor $P_\infty$ ensures that $P_\infty E_0$ and $\|\nabla u_0\|_{L^2}^2$ have the same scaling, while the prefactor $\rho_\infty^{1-\gamma}$ makes their sum invariant under
\eqref{2sca}.
\vspace{2mm}

We state our first result.
\begin{theorem}[Far-field non-vacuum]\label{th1}
Let $\rho_\infty>0$. Assume that the initial data $(\rho_0\geq 0,u_0)$ satisfies
\begin{align}\label{d1}
G(\rho_{0}\mid\rho_\infty)\in L^{1}, \quad  (\rho_{0} -\rho_{\infty})\in H^{1} \cap W^{1,q}, \quad \nabla u_{0} \in L^{2}, \quad \sqrt{\rho_{0}} u_{0} \in L^{2}.
\end{align}
There exists $\varepsilon_{\mathrm nv}>0$ depending only on $A,\gamma,\mu,\lambda$, and $\mathcal{R}$ such that if
\begin{equation}\label{S0}
\mathcal S_0^{\mathrm{nvac}}\le\varepsilon_{\mathrm nv},
\end{equation}
the Cauchy problem \eqref{1}--\eqref{3} possesses a unique global strong solution.
\end{theorem}

We now consider the far-field vacuum case $\rho_\infty=0$. Recall the intrinsic length scale $L_0$ and the normalized moment $\mathcal M_0$ given by
\begin{align*}
L_0=\left(\frac{m_0}{\bar\rho}\right)^{1/2},\quad\mathcal M_0=\frac{1}{m_0L_0}
\int\rho_0(x)\sqrt{L_0^2+|x|^2}\,dx.
\end{align*}
In this case, the initial energy takes the form
\begin{align*}
E_0=\int\bigg(\frac12\rho_0|u_0|^2+\frac{A}{\gamma-1}\rho_0^\gamma\bigg)dx,
\end{align*}
and we introduce the scaling-invariant initial quantity
\begin{align}\label{SS}
\mathcal S_0^{\mathrm{vac}}\triangleq\mathcal M_0^{1/3}\big[1+(Am_0\bar\rho^\gamma)^{-2/3}\big]\bar\rho
\left[E_0+m_0\left(\|\nabla u_0\|_{L^2}^2+\|P(\rho_0)\|_{L^2}^2\right)\right].
\end{align}

We are in a position to state our second main result.
\begin{theorem}[Far-field vacuum]\label{th2}
Let $\rho_\infty=0$. Assume that the initial data $(\rho_{0}\ge 0,u_{0})$ satisfies
\begin{align}\label{ini}
\bar{x}^{a} \rho_{0} \in L^{1} \cap H^{1} \cap W^{1, q}, \quad \nabla u_{0} \in L^{2}, \quad \sqrt{\rho} u_{0} \in L^{2}.
\end{align}
There exists $\varepsilon_{\mathrm v}>0$ depending only on $A,\gamma,\mu$, and $\lambda$ such that if
\begin{align}\label{S1}
\mathcal S_0^{\mathrm{vac}}\le\varepsilon_{\mathrm v},
\end{align}
the Cauchy problem \eqref{1}--\eqref{3} admits a unique global strong solution.
\end{theorem}

Several remarks are in order.

\begin{remark}
The quantities defined in \eqref{SS1} and \eqref{SS} are invariant under the scaling transformation \eqref{2sca}; see subsection \ref{dd} for details.
\end{remark}

\begin{remark}[Related independent results]
Our Theorems \ref{th1} and \ref{th2} presented in this work are the results announced in \cite[Remark 1.7]{XuZhong202611}. During the completion of the present work, we find that Wen--Wu \cite{WW2026} and L\"u--Song \cite{LS2026} independently established the global existence and uniqueness of strong solutions for \eqref{1} with far-field vacuum under different scaling-invariant smallness conditions. It is worth mentioning that the scaling-invariant initial quantities introduced in these two works differ substantially in form from the quantity appearing in \eqref{S1}. In particular, our condition is formulated in terms of the initial energy, first-order
dissipation, total mass, and a quantity associated with the first spatial moment of the density.
Furthermore, unlike their results focusing on the far-field vacuum case, the present paper treats both far-field vacuum and non-vacuum states within a unified scaling-invariant framework.
\end{remark}

\begin{remark}[Different treatments of the pressure]
An important difference between the far-field vacuum and non-vacuum problems lies in the treatment
of the pressure terms in the a \textit{priori} estimates.
	
For the non-vacuum case $\rho_\infty>0$. Under the density assumption $0\leq\rho\leq \bar\rho$, we have, for $\RR=\bar\rho/\rho_\infty$,
\begin{align*}
\|P(\rho)\|_{L^\infty}\leq A(\RR)^\gamma\rho_\infty^\gamma,\quad
\|P(\rho)-P_\infty\|_{L^\infty}\leq C(\RR,A,\gamma)\rho_\infty^\gamma.
\end{align*}
Moreover, the relative potential energy provides the estimate, see \eqref{eqeq2},
\begin{align*}
\|P(\rho)-P_\infty\|_{L^2}^2\leq C(\RR,A,\gamma)\rho_\infty^\gamma\int G(\rho\mid\rho_\infty)\,dx.
\end{align*}
Thus, the pressure contributions can be controlled in terms of $\rho_\infty$ and the relative
energy. In particular, the BMO-based interpolation estimate \eqref{L4} allows the contribution of the background pressure to retain the dissipative factor $\|\nabla u\|_{L^2}^2$. This structure is essential for closing
the lower-order and time-weighted estimates \eqref{eqlow} and \eqref{w0}.
	
The situation is different when $\rho_\infty=0$. In this case, the pressure  equation
\begin{align*}
P_t+u\cdot\nabla P+\gamma P\operatorname{div}u=0
\end{align*}
provides an additional dissipative mechanism. Indeed, testing the above equation by suitable powers of $P$, together with the effective viscous flux decomposition, yields the cubic and quartic pressure estimates \eqref{eq42} and \eqref{65}. These estimates play crucial roles in controlling the nonlinear terms and establishing the time-weighted bounds \eqref{60} and \eqref{61}.
	
Consequently, although the two far-field problems share the same  scaling transformation \eqref{2sca}, their pressure estimates are organized differently. The non-vacuum argument exploits the positive background
state and relative energy, whereas the vacuum argument makes essential use of the dissipative structure of the pressure evolution equation.
\end{remark}

\begin{remark}[The role of the far-field density and total mass]
For the non-vacuum case $\rho_\infty>0$, the relative energy and the weighted Poincar\'e inequality \eqref{eqw} allow us to recover unweighted $L^2$ control from density-weighted quantities. The ratio $\RR=\bar\rho/\rho_\infty$ is preserved under the scaling \eqref{2sca}, and the powers of $\rho_\infty$ appearing in the estimates are retained explicitly. In particular, the intrinsic time weight
$\min\{1,\rho_\infty^\gamma t\}$ in \eqref{w0} ensures that the corresponding estimates remain compatible with the scaling.
	
For the vacuum case $\rho_\infty=0$, this mechanism is no longer available. Instead, the finite total mass
and the first spatial moment provide the necessary localization of mass. Combining this property with a weighted Sobolev inequality yields the density-weighted $L^p$ estimate in Lemma \ref{zuilem}. Furthermore, the explicit mass factor in the pressure estimate \eqref{eqp1} determines the appearance of $m_0D_{0}$ in the initial smallness quantity. The mass also provides the scale-adapted time parameter used in the weights $m_0+t$ and $t(m_0+t)$; see \eqref{60} and \eqref{61}.
	
These quantities provide two different mechanisms for overcoming the low-frequency difficulty of the two-dimensional Cauchy problem. More importantly, their scaling properties determine the appropriate combinations of energy, dissipation, and time weights, leading to the scaling-invariant initial conditions in the two far-field settings.
\end{remark}

\begin{remark}[Weighted assumptions on the initial density]
In the far-field vacuum case, the global a priori estimates in Section \ref{sec4} use the finite first spatial moment
\begin{align}\label{nn}
\int \rho_0(x)|x|\,dx<\infty
\end{align}
as an additional spatial decay assumption on the initial density, see \eqref{kk}. The local existence theory in \cite{LiLiang2014} assumes the weighted regularity condition $\bar{x}^{a}\rho_0\in L^1\cap H^1\cap W^{1,q}$, which implies \eqref{nn}. In our argument, the global density bootstrap uses the first-moment condition as an additional assumption on the spatial decay of the initial density. Thus, it is natural to ask whether the local existence theory can be established under weaker spatial decay and weighted regularity assumptions. We leave this question for future investigation.
\end{remark}


\begin{remark}[Two distinct far-field mechanisms]
Although \eqref{S0} and \eqref{S1} have different forms, they follow from the same scaling \eqref{2sca}. The vacuum criterion involves finite mass and the first moment, which supply the localization information needed for the low-frequency analysis. The non-vacuum criterion uses the relative potential energy and the background pressure. In particular, the vacuum result is not obtained by simply setting $\rho_\infty=0$ in the non-vacuum condition; the change of far-field state fundamentally changes the available mass and structures.
\end{remark}

\begin{remark}[Vacuum and large pointwise oscillations]
In both theorems, the initial density $\rho_0$ is allowed to vanish at an interior region. Moreover, the smallness condition \eqref{S0} imposes no smallness requirement on $\|\rho_0-\rho_\infty\|_{L^\infty}$. The admissibility of vacuum and large oscillations is compatible with the earlier small-energy theory \cite{HLX2012,LiXin2019}. The present emphasis is the scaling-invariant formulation and the corresponding estimates rather than the mere presence of vacuum.
\end{remark}

\begin{remark}[A possible extension to large initial energy]
An interesting direction for further investigation is to extend the present results to the regime of large initial energy. In this connection, Hong, Hou, Peng, and Zhu \cite{HHPZ2024} established the global existence of classical solutions with large initial energy and vacuum to the three-dimensional compressible Navier--Stokes equations when the adiabatic exponent is sufficiently close to $1$, allowing either far-field vacuum or a sufficiently small positive far-field density.
	
In the present two-dimensional setting, the smallness thresholds in both far-field cases depend on $\gamma$. This naturally raises the question of whether the a priori estimates can be refined to obtain smallness thresholds independent of $\gamma$ near $1$. More importantly, it would be desirable to identify a refined scaling-invariant initial condition that exploits the smallness of $\gamma-1$ and allows the initial energy to be large. Establishing such Nishida--Smoller-type global solutions in two dimensions, for both far-field vacuum and non-vacuum states, is an interesting problem that we intend to investigate in future work.	
\end{remark}

\subsection{Difficulties and strategy of the proof}

We briefly explain the main difficulties and key ideas in the proofs of our main results. The first difficulty arises from the lack of a suitable
homogeneous Sobolev embedding in two dimensions. In three dimensions, the classical embedding $\dot H^1(\mathbb R^3)\hookrightarrow L^6(\mathbb R^3)$
provides the estimate $\|u\|_{L^6(\mathbb R^3)}\le C\|\nabla u\|_{L^2(\mathbb R^3)}$, which plays an important role in controlling nonlinear terms involving the velocity and momentum. In contrast, in two dimensions, $\|\nabla u\|_{L^2}$ alone does not control $\|u\|_{L^p}$ for any finite $p$ without additional low-frequency
information. Moreover, the basic energy inequality controls $\sqrt{\rho}u$ rather than $u$ itself, and this distinction is particularly significant
when the initial density is allowed to vanish. Consequently, the estimates of the velocity, momentum, and effective viscous flux used in three dimensions cannot be directly carried over to the whole plane. The second difficulty is to obtain a uniform-in-time upper bound for the density under a smallness
condition that is invariant under the exact scaling \eqref{2sca}. In particular, the initial energy and the first-order dissipation have different scaling orders, and the constants in the nonlinear estimates must be tracked carefully. Furthermore, the far-field vacuum and non-vacuum problems have
different low-frequency structures, so the same energy functional and velocity estimates cannot be used directly in both cases.

To overcome these difficulties, we combine scaling-invariant energy estimates with a bootstrap argument for the density. Along particle trajectories, the continuity equation can be written as
\begin{align*}
\frac{d}{dt}\rho=-\frac{\rho}{2\mu+\lambda}\bigl(F+P(\rho)-P_\infty\bigr).
\end{align*}
Thus, suitable time-integrated estimates for $F$, together with the damping effect of the pressure, allow us to apply Zlotnik's inequality (see Lemma \ref{zo}) and to derive the required density upper bound. Therefore, the main issue is to establish sufficiently strong estimates for the effective viscous flux, both near the initial time and over arbitrarily long time intervals.

For the \emph{far-field vacuum case} $\rho_\infty=0$, the principal obstacle is the lack of unweighted $L^p$ estimates for the velocity. To this end, we exploit the finite total mass and the first spatial moment of the initial density to obtain quantitative mass localization; see \eqref{70} and \eqref{71}. Combining this property with a weighted Sobolev inequality yields the density-weighted $L^p$ estimate in Lemma \ref{zuilem}. This provides the necessary control of $\rho\dot u$ in higher Lebesgue spaces without requiring an unweighted $L^2$ bound for the velocity. Another important ingredient is the pressure estimate \eqref{eqp1}, which contains the explicit factor $m_0^{1/2}$. Together with \eqref{eqp2}, this estimate naturally leads to the scaling-invariant combination of the initial energy and first-order dissipation, see Remark \ref{rem}. We then establish the time-weighted estimates \eqref{60} and \eqref{61}, with the intrinsic weights $m_0+t$ and $t(m_0+t)$. These estimates provide algebraic decay of the velocity gradient and material derivative, which is essential for controlling the effective viscous flux over large time intervals. More precisely, the interpolation arguments in Lemma \ref{flux} yield the short-time $L^1_tL^\infty_x$ estimate \eqref{90} and the large-time $L^{3/2}_tL^\infty_x$ estimate \eqref{91}. The latter estimate, combined with the pressure damping in the density equation, is sufficient to apply Zlotnik's inequality and improve the density bootstrap; see Proposition \ref{prop}. Tracking the dependence on the initial mass, the density upper bound, and the first spatial moment leads to the refined scaling-invariant smallness condition \eqref{bao}.

For the \emph{far-field non-vacuum case} $\rho_\infty>0$, the positive background density provides a different mechanism for recovering low-frequency control. The relative potential energy yields the estimates \eqref{eqeq1} and \eqref{eqeq2}, while the weighted Poincar\'e inequality \eqref{eqw} allows us to control the unweighted $L^2$ norm of the material derivative in terms of density-weighted quantities and its spatial gradient. A further difficulty is caused by the nonzero background pressure $P_{\infty}$, whose contribution to the nonlinear estimates cannot be controlled by a direct time integration without additional decay information. To overcome this obstacle, we exploit the two-dimensional embedding $\dot H^1(\mathbb R^2)\hookrightarrow
\mathrm{BMO}(\mathbb R^2)$, together with elliptic estimates for the effective viscous flux and vorticity. The resulting interpolation estimate \eqref{L4} takes the form
\begin{align*}
\|\nabla u\|_{L^4}^4\le C\|\nabla u\|_{L^2}^2
\left(\rho_\infty\|\sqrt\rho\dot u\|_{L^2}^2+\rho_\infty^{2\gamma}\right).
\end{align*}
An important feature of this inequality is that the background-pressure contribution remains multiplied by the dissipative factor $\|\nabla u\|_{L^2}^2$, whose time integral is controlled by the basic energy inequality; see Remark \ref{u1}. Based on this estimate, we derive the lower-order bound \eqref{eqlow} and then establish the time-weighted material-derivative estimate \eqref{w0} via
\begin{align*}
\sigma(t)=\min\{1,\rho_\infty^\gamma t\}.
\end{align*}
This choice is compatible with the exact scaling and ensures that the
constants in the weighted estimates are independent of the absolute value
of $\rho_\infty$ (see Remark \ref{weight} for details). The resulting estimates yield short-time $L^1_tL^\infty_x$ control of $F$ in \eqref{eqsh} and large-time $L^3_tL^\infty_x$ control in \eqref{eqsh1}. These bounds allow us to control the increments of the effective viscous flux along particle trajectories and close the density bootstrap through Zlotnik's inequality;
see Proposition \ref{boot}. Unlike the far-field vacuum argument, this approach does not require a separate algebraic decay estimate for the velocity gradient.

In summary, the two far-field problems share the same underlying strategy:
scaling-invariant energy estimates, material-derivative bounds, effective
viscous flux estimates, and a density bootstrap based on Zlotnik's inequality.
Their main difference lies in the treatment of low frequencies and the resulting time-integrability properties of the effective viscous flux. The far-field vacuum case relies on finite-mass localization and weighted velocity estimates, whereas the non-vacuum case uses relative energy and BMO interpolation. These two mechanisms provide a unified framework for obtaining the uniform density bounds required in the global strong-solution analysis.

The remainder of the paper is organized as follows. Section \ref{sec2} collects auxiliary inequalities and establishes scaling exponents and scaling-invariant initial quantities. Section \ref{sec3} is devoted to dealing with the far-field non-vacuum case, while we handle the far-field vacuum case in Section \ref{sec4}.

\section{Preliminaries}\label{sec2}

In this section we collect auxiliary results and provide scaling exponents and scaling-invariant initial quantities.

\subsection{Some inequalities and blow-up criterion}
\begin{lemma}[Gagliardo--Nirenberg inequality]\label{gn}
For $s \in[2, \infty)$, $q \in(1, \infty)$, and $r \in(2, \infty)$, there exists some generic constant $C>0$  which may depend on $s, q$, and $r$ such that, for $f \in H^{1}(\mathbb{R}^{2})$ and $g \in L^{q}\left(\mathbb{R}^{2}\right) \cap D^{1, r}(\mathbb{R}^{2})$,
\begin{gather}
\|f\|_{L^{s}} \leq C\|f\|_{L^{2}}^{2 / s}\|\nabla f\|_{L^{2}}^{(s-2) / s}, \label{g1}\\[1mm]
\|g\|_{C(\overline{\mathbb{R}}^{2})} \leq C\|g\|_{L^{q}}^{q(r-2) /(2 r+q(r-2))}\|\nabla g\|_{L^{r}}^{2 r /(2 r+q(r-2))}.\label{g2}
\end{gather}
\end{lemma}
\begin{proof}
See \cite[Theorem 2.2]{L1968}.
\end{proof}
	
\begin{lemma}[Generalised Gagliardo--Nirenberg inequality]\label{ggn}
For any  $1<q<r<\infty$, if $f \in L^{q, \infty}(\mathbb{R}^{2}) \cap \operatorname{BMO}(\mathbb{R}^{2})$, then  $f \in L^{r}(\mathbb{R}^{2})$ and there exists a constant $C=C(q, r)$ such that
\begin{align*}
\|f\|_{L^{r}} \leq C\|f\|_{L^{q, \infty}}^{q / r}[f]_{\mathrm{BMO}}^{1-q / r}.
\end{align*}
\end{lemma}
\begin{proof}
See \cite[Theorem 9.1]{MRR2013}.
\end{proof}	
	
\begin{lemma}[Zlotnik's inequality]\label{zo}
Let the function $y$ satisfy
\begin{align*}
y^{\prime}(t)=g(y)+b^{\prime}(t) \text { on }[0, T],\quad  y(0)=y^{0},
\end{align*}
with $g \in C(R)$ and $y$, $b \in W^{1,1}(0, T)$. If $g(\infty)=-\infty$ and
\begin{align*}
b\left(t_{2}\right)-b\left(t_{1}\right) \leq N_{0}+N_{1}\left(t_{2}-t_{1}\right),
\end{align*}
for all $0 \leq t_{1}<t_{2} \leq T$  with some $N_{0} \geq 0$ and $N_{1} \geq 0$, then
\begin{align}\label{zo1}
y(t) \leq \max \left\{y^{0}, \bar{\zeta}\right\}+N_{0}<\infty\quad \text { on }\quad [0, T],
\end{align}
where $\bar{\zeta}$ is a constant such that
\begin{align*}
g(\zeta) \leq-N_{1} \quad\text { for }\quad \zeta \geq \bar{\zeta}.
\end{align*}
\end{lemma}
\begin{proof}
See \cite[Lemma 1.3]{Z2000}.
\end{proof}

\begin{lemma}[Hardy--Littlewood--Sobolev inequality]\label{hls}
Let $q, r>1$ and $0<\lambda<n$ with $1 / q+\lambda / n+1 / r=2$. Let $f \in L^{q}\left(\mathbb{R}^{n}\right)$  and  $h \in L^{r}\left(\mathbb{R}^{n}\right)$. Then there exists a sharp constant $C(n, \lambda, q)$ independent of $f$ and $h$ such that
\begin{align}\label{ff}
\left|\int_{\mathbb{R}^{n}} \int_{\mathbb{R}^{n}} f(x) | x- y |^{-\lambda} h(y) d x  d y \right|  \leq C(n, \lambda, q)\|f\|_{q}\|h\|_{r}.
\end{align}
\end{lemma}
\begin{proof}
See \cite[Theorem 4.3]{LiebLoss2001}.
\end{proof}
	
\begin{remark}\label{br}
In \eqref{ff}, taking $n=2, \lambda=1, q=\frac{4\gamma}{4\gamma-1}$, and $r=\frac{4\gamma}{2\gamma+1}$, we have
\begin{align*}
\frac1q+\frac12+\frac1r=2,\quad 1<r<2.
\end{align*}
Hence, by $L^{4\gamma}$ duality, one gets that
\begin{align}\label{eqHLS}
\left\|\int\frac{h(y)}{|x-y|}\,dy\right\|_{L^{4\gamma}}
\le C\|h\|_{L^{\frac{4\gamma}{2\gamma+1}}}.
\end{align}
Furthermore, the fundamental solution of
$-\Delta$ in $\mathbb R^2$ is
$\Gamma(x)=-(2\pi)^{-1}\log|x|$.
Consequently,
\begin{align*}
\bigl|(-\Delta)^{-1}\operatorname{div}h(x)\bigr|=\left|\int\nabla\Gamma(x-y)\cdot h(y)\,dy\right| \le C\int
\frac{|h(y)|}{|x-y|}\,dy.
\end{align*}
Applying \eqref{eqHLS} with $h=\rho\dot u$ and then using H\"older's inequality, we conclude that
\begin{align*}
\bigl\|(-\Delta)^{-1}\operatorname{div}(\rho\dot u)\bigr\|_{L^{4\gamma}}
\le C\|\rho\dot u\|_{L^{\frac{4\gamma}{2\gamma+1}}}
\le C\|\rho\|_{L^{2\gamma}}^{1/2}\|\sqrt\rho \dot u\|_{L^{2}}.
\end{align*}
\end{remark}	

\begin{lemma}[Weighted inequality]\label{lious}
For  $m \in[2, \infty)$  and  $\theta \in(1+m / 2, \infty)$, there exists a positive constant  $C$  such that for all  $v \in D^{1}\left(\mathbb{R}^{2}\right)$,
\begin{align*}
\left(\int_{\mathbb{R}^{2}} \frac{|v|^{m}}{e+|x|^{2}}\left(\log \left(e+|x|^{2}\right)\right)^{-\theta}\, d x\right)^{1 / m} \leq C\|v\|_{L^{2}\left(B_{1}\right)}+C\|\nabla v\|_{L^{2}\left(\mathbb{R}^{2}\right)}.
\end{align*}	
\end{lemma}
\begin{proof}
See \cite[Theorem B.1]{P1996}.
\end{proof}
	
\begin{lemma}[Blow-up criterion]\label{blowup}
Let $(\rho,u)$ be a strong solution to the Cauchy problem \eqref{1}--\eqref{3}, with the initial data satisfying
\eqref{d1} or \eqref{ini}, according to whether $\rho_\infty>0$ or $\rho_\infty=0$. Let $T^*>0$ be the maximal existence time of the strong solution. If $T^*<\infty$, then
\begin{align}\label{eq:blowup}
\lim_{T\rightarrow T^*}\|\rho\|_{L^\infty(0,T;L^\infty)}=\infty.
\end{align}
\end{lemma}	
\begin{proof}
For the case $\rho_\infty=0$, it follows from \cite[Theorem 1.1]{J2019} that there exists a constant $s_0>1$ such that
\begin{align*}
\lim_{T\rightarrow T^*}\|\rho\|_{L^\infty(0,T;L^{s_0})}=\infty.
\end{align*}
By conservation of mass, for each $t\in[0,T^*)$ we have
\begin{align*}
\|\rho\|_{L^{s_0}}\leq m_0^{1/s_0}\|\rho\|_{L^\infty}^{1-1/s_0},
\end{align*}
which immediately yields \eqref{eq:blowup}.
		
When $\rho_\infty>0$, the desired result follows by a slight modification of the proof of \cite[Theorem 1.2]{W2026}.
\end{proof}

\subsection{Determination of the scaling exponents}\label{scax}

For any $\ell>0$, consider the transformation
\begin{align*}
\rho^\ell(x,t)\triangleq\ell^\alpha\rho(\ell x,\ell^\kappa t),\quad
u^\ell(x,t)\triangleq\ell^\beta u(\ell x,\ell^\kappa t),
\end{align*}
where the exponents $\alpha$, $\beta$, and $\kappa$ are to be determined.

By the chain rule, the two terms in the continuity equation \eqref{1}$_1$ scale as
\begin{align*}
\partial_t\rho^\ell=\ell^{\alpha+\kappa}\rho_t(\ell x,\ell^\kappa t),\quad\operatorname{div}(\rho^\ell u^\ell)
=\ell^{\alpha+\beta+1}\operatorname{div}(\rho u)(\ell x,\ell^\kappa t),
\end{align*}
and the terms in the momentum equation \eqref{1}$_2$ transform as
\begin{gather*}
\rho^\ell\partial_tu^\ell=\ell^{\alpha+\beta+\kappa}(\rho u_t)(\ell x,\ell^\kappa t),\quad
\rho^\ell(u^\ell\cdot\nabla)u^\ell=\ell^{\alpha+2\beta+1}\rho(u\cdot\nabla)u(\ell x,\ell^\kappa t),\\
\Delta u^\ell=\ell^{\beta+2}\Delta u(\ell x,\ell^\kappa t),\quad\nabla\operatorname{div}u^\ell
=\ell^{\beta+2}\nabla\operatorname{div}u(\ell x,\ell^\kappa t),\quad\nabla P(\rho^\ell)
=\ell^{\gamma\alpha+1}\nabla P(\rho)(\ell x,\ell^\kappa t).
\end{gather*}
Thus, owing to invariance of the continuity equation and the momentum equation, one needs to have
\begin{align*}
\begin{cases}
\kappa=\beta+1,\\
\alpha+\beta+\kappa=\alpha+2\beta+1=\beta+2=\gamma\alpha+1,
\end{cases}
\end{align*}
which yields
\begin{align*}
\alpha=\frac{2}{\gamma+1},\quad \beta=\frac{\gamma-1}{\gamma+1},\quad \kappa=\frac{2\gamma}{\gamma+1}.
\end{align*}
Consequently, the scaling transformation is
\begin{align}\label{sc}
\rho^\ell(x,t)=\ell^{\frac{2}{\gamma+1}}\rho\big(\ell x,\ell^{\frac{2\gamma}{\gamma+1}}t\big),
\quad u^\ell(x,t)=\ell^{\frac{\gamma-1}{\gamma+1}}u\big(\ell x,\ell^{\frac{2\gamma}{\gamma+1}}t\big).
\end{align}

\subsection{Scaling properties of the initial energy and dissipation}

Under \eqref{sc}, the far-field density and pressure transform have the forms
\begin{align}\label{s1}
\rho_\infty^\ell&=\ell^{\frac{2}{\gamma+1}}\rho_\infty,
\quad P_\infty^\ell=\ell^{\frac{2\gamma}{\gamma+1}}P_\infty.
\end{align}
The corresponding initial density upper bound satisfies
\begin{align}\label{s2}
\bar\rho^\ell=\ell^{\frac{2}{\gamma+1}}\bar\rho,
\quad\frac{\bar\rho^\ell}{\rho_\infty^\ell}=\frac{\bar\rho}{\rho_\infty}=\mathcal{R}\quad\text{if }\rho_\infty>0.
\end{align}

Recall that
\begin{align*}
D_{0}=\|\nabla u_0\|_{L^2}^2+\|P(\rho_0)-P_\infty\|_{L^2}^2,\quad
E_0=\int\bigg(\frac12\rho_0|u_0|^2+G(\rho_0\mid\rho_\infty)\bigg)dx,
\end{align*}
where
\begin{align*}
G(\rho\mid\rho_\infty)=\frac{A}{\gamma-1}\left[\rho^\gamma-\gamma\rho_\infty^{\gamma-1}\rho+(\gamma-1)\rho_\infty^\gamma\right].
\end{align*}
Since the spatial dimension is two, the change of variables $y=\ell x$ gives $dx=\ell^{-2}dy$, which combined with \eqref{sc} indicates
\begin{gather}\label{s5}
\|\nabla u_0^\ell\|_{L^2}^2=\ell^{\frac{2(\gamma-1)}{\gamma+1}}\|\nabla u_0\|_{L^2}^2,\quad\|P(\rho_0^\ell)-P_\infty^\ell\|_{L^2}^2
=\ell^{\frac{2(\gamma-1)}{\gamma+1}}\|P(\rho_0)-P_\infty\|_{L^2}^2,
\\ \label{s3}
\int\frac12\rho_0^\ell|u_0^\ell|^2\,dx=\ell^{-\frac{2}{\gamma+1}}\int\frac12\rho_0|u_0|^2\,dx.
\end{gather}

Moreover, the relative potential energy density is homogeneous of degree $\gamma$ in
$(\rho,\rho_\infty)$, namely,
\begin{align*}
G(h\rho\mid h\rho_\infty)=h^\gamma G(\rho\mid\rho_\infty),\quad h>0.
\end{align*}
Hence, taking $h=\ell^{2/(\gamma+1)}$ and changing variables, we obtain
\begin{align*}
\int G(\rho_0^\ell\mid\rho_\infty^\ell)\,dx=\ell^{-\frac{2}{\gamma+1}}\int G(\rho_0\mid\rho_\infty)dx,
\end{align*}
which together with \eqref{s5} and \eqref{s3} leads to
\begin{align}\label{s6}
E_0^\ell=\ell^{-\frac{2}{\gamma+1}}E_0,\quad D_{0}^\ell=\ell^{\frac{2(\gamma-1)}{\gamma+1}}D_{0}.
\end{align}

\subsection{Scaling invariance}\label{dd}

Suppose $\rho_\infty>0$. Owing to $\mathcal S_0^{\mathrm{nvac}}=\rho_\infty^{1-\gamma}\left[D_{0}+P_\infty E_0\right]$,
one obtains from \eqref{s1} and \eqref{s6} that
\begin{align*}
D_{0}^\ell +P_\infty^\ell E_0^\ell=\ell^{\frac{2(\gamma-1)}{\gamma+1}}\left[D_{0}+P_\infty E_0\right].
\end{align*}
Hence
\begin{align*}
\left(\mathcal S_0^{\mathrm{nvac}}\right) ^{\ell}=(\rho_\infty^\ell)^{1-\gamma}\big[D^\ell_{0}+P_\infty^\ell E_0^\ell\big]
=\ell^{\frac{2(1-\gamma)}{\gamma+1}}\ell^{\frac{2(\gamma-1)}{\gamma+1}}\mathcal S_0^{\mathrm{nvac}}=\mathcal S_0^{\mathrm{nvac}}.
\end{align*}
Since $\mathcal{R}^\ell=\mathcal{R}$ by \eqref{s2}, the smallness condition
\begin{align*}
\mathcal S_0^{\mathrm{nvac}}\leq\varepsilon_*(\mathcal{R},A,\gamma,\mu,\lambda)
\end{align*}
is invariant under \eqref{sc}.

Now suppose $\rho_\infty=0$. Recall that
\begin{align*}
m_0=\int \rho_0\,dx,\quad \Pi_0=Am_0\bar\rho^\gamma,\quad
L_0=\Big(\frac{m_0}{\bar\rho}\Big)^{1/2},\quad
\mathcal M_0=\frac{1}{m_0L_0}\int\rho_0(x)\sqrt{L_0^2+|x|^2}\,dx.
\end{align*}
Define $m_0^\ell$, $\Pi_0^\ell$, $L_0^\ell$, and $\mathcal M_0^\ell$ using the scaled initial
density and its upper bound. By \eqref{sc} and \eqref{s2}, we have
\begin{align}\label{mo}
m_0^\ell=\ell^{-\frac{2\gamma}{\gamma+1}}m_0,\quad \Pi_0^\ell=A m_0^\ell(\bar\rho^\ell)^\gamma=A\ell^{-\frac{2\gamma}{\gamma+1}}
m_0\ell^{\frac{2\gamma}{\gamma+1}}\bar\rho^\gamma=\Pi_0,\quad L_0^\ell=\ell^{-1}L_0.
\end{align}
Moreover, the change of variables gives
\begin{align*}
\int \rho_{0}^{\ell}(x) \sqrt{\left(L_{0}^{\ell}\right)^{2}+|x|^{2}}\, dx
&=\ell^{\frac{2}{\gamma+1}-2} \int  \rho_{0}(y) \sqrt{\ell^{-2} L_{0}^{2}+\ell^{-2}|y|^{2}}\, dy \\
&=\ell^{\frac{2}{\gamma+1}-3} \int  \rho_{0}(y) \sqrt{L_{0}^{2}+|y|^{2}}\, dy \\
&=\ell^{-\frac{3 \gamma+1}{\gamma+1}} \int  \rho_{0}(y) \sqrt{L_{0}^{2}+|y|^{2}}\, dy.
\end{align*}
Since $m_0^\ell L_0^\ell=\ell^{-\frac{3\gamma+1}{\gamma+1}}m_0L_0$, we conclude that $\mathcal M_0^\ell=\mathcal M_0$.

Furthermore, \eqref{s2}, \eqref{s6}, and \eqref{mo} show that both $\bar\rho E_0$ and $\bar\rho m_0D_0$ are invariant under \eqref{sc}, which implies the quantity
$\mathcal S_0=\bar\rho\left( E_0+m_0D_0\right)$
satisfies $\mathcal S_0^\ell=\mathcal S_0$. Consequently, the initial quantity
\begin{align*}
\mathcal S_0^{\mathrm{vac}}=\mathcal M_0^{1/3}\big(1+\Pi_0^{-2/3}\big)\mathcal S_0
\end{align*}
is also invariant. This completes the verification of the scaling invariance of the initial quantities in both far-field regimes.

\section{Compressible Navier--Stokes equations with far-field non-vacuum}\label{sec3}
For $\rho_\infty>0$ and initial data $(\rho_0,u_0)$ satisfying \eqref{d1}, it follows from \cite{Luo2012} that the problem \eqref{1}--\eqref{3} admits a unique local strong solution on $\mathbb{R}^2\times(0,T]$ for some $T>0$. In this section we establish the {\it a priori} estimates for such solution. Throughout this section, $c$ and $C$ denote the generic positive constant depending only on $\RR$, $A$, $\gamma$, $\mu$, and $\lambda$, unless
otherwise specified. Their values may vary from line to line.

\subsection{Basic estimates}\label{subsec20}
\begin{lemma}\label{2lem1}
Let $v\in H^1$ and $(\rho-\rho_\infty)\in L^2$, then
\begin{align}\label{eqw}
\|v\|_{L^2}\leq C\rho_\infty^{-1/2}\|\sqrt\rho v\|_{L^2}+C\rho_\infty^{-1}
\|\rho-\rho_\infty\|_{L^2}\|\nabla v\|_{L^2}.
\end{align}
In particular, if $\dot u\in H^1$, then
\begin{gather}
\|\dot u\|_{L^2}\leq C\rho_\infty^{-1/2}\|\sqrt\rho \dot u\|_{L^2}+C\rho_\infty^{-1}\|\rho-\rho_\infty\|_{L^2}
\|\nabla\dot u\|_{L^2},\label{eqd1}\\
\|\dot u\|_{L^4}\leq C\rho_\infty^{-1/4}\|\sqrt\rho \dot u\|_{L^2}^{1/2}
\|\nabla\dot u\|_{L^2}^{1/2}+C\rho_\infty^{-1/2}\|\rho-\rho_\infty\|_{L^2}^{1/2}\|\nabla\dot u\|_{L^2}.
\label{eqd2}
\end{gather}
\end{lemma}
\begin{proof}
It follows from H\"older's inequality and Gagliardo--Nirenberg inequality that
\begin{align}\label{xz1}
\int \rho|v|^2\,dx\leq \rho_\infty\|v\|_{L^2}^2+\|\rho-\rho_\infty\|_{L^2}
\|v\|_{L^4}^2\leq \rho_\infty\|v\|_{L^2}^2+C\|\rho-\rho_\infty\|_{L^2}\|v\|_{L^2}\|\nabla v\|_{L^2}.
\end{align}
Writing $\rho_\infty=\rho+(\rho_\infty-\rho)$, we obtain from H\"older's inequality, Gagliardo--Nirenberg inequality, and Young's inequality that
\begin{align*}
\rho_\infty\|v\|_{L^2}^2&=\|\sqrt\rho v\|_{L^2}^2+\int(\rho_\infty-\rho)|v|^2\,dx\\
&\leq \|\sqrt\rho v\|_{L^2}^2+\|\rho-\rho_\infty\|_{L^2}\|v\|_{L^4}^2\\
&\leq \|\sqrt\rho v\|_{L^2}^2+C\|\rho-\rho_\infty\|_{L^2}\|v\|_{L^2}\|\nabla v\|_{L^2}\\
&\leq \|\sqrt\rho v\|_{L^2}^2+\frac{\rho_\infty}{2}\|v\|_{L^2}^2
+C\rho_\infty^{-1}\|\rho-\rho_\infty\|_{L^2}^2\|\nabla v\|_{L^2}^2,
\end{align*}
which combined with \eqref{xz1} implies \eqref{eqw}.
	
Taking $v=\dot u$ in \eqref{eqw} immediately yields \eqref{eqd1}. Furthermore, applying Gagliardo--Nirenberg inequality and \eqref{eqd1}, we find
\begin{align*}
\|\dot u\|_{L^4} \leq C\|\dot u\|_{L^2}^{1/2}\|\nabla\dot u\|_{L^2}^{1/2} \leq C\rho_\infty^{-1/4}
\|\sqrt\rho\,\dot u\|_{L^2}^{1/2}\|\nabla\dot u\|_{L^2}^{1/2} +C\rho_\infty^{-1/2}
\|\rho-\rho_\infty\|_{L^2}^{1/2}\|\nabla\dot u\|_{L^2},
\end{align*}
proving \eqref{eqd2}.
\end{proof}

We recall the scaling-invariant quantity
\begin{align*}
\mathcal S_0^{\mathrm{nvac}}=\rho_\infty^{1-\gamma}\left(\|\nabla u_0\|_{L^2}^2+P_\infty E_{0}\right).
\end{align*}

\begin{lemma}\label{lem1}
It holds that
\begin{align}\label{eqba}
\sup_{0\leq t\leq T}E(t)+\mu\int_0^T\|\nabla u\|_{L^2}^2\, dt
\leq C\mathcal S_0^{\mathrm{nvac}}\rho_\infty^{-1}.	
\end{align}
If, in addition, $0\leq\rho(x,t)\leq\bar\rho=\RR\rho_\infty$, then
\begin{gather}
c\rho_\infty^{\gamma-2}(\rho-\rho_\infty)^2\leq G(\rho\mid\rho_\infty)
\leq C\rho_\infty^{\gamma-2}(\rho-\rho_\infty)^2,\label{eqeq1}\\
p^2\leq C\rho_\infty^\gamma G(\rho\mid\rho_\infty),\label{eqeq2}\\
\sup_{0\leq t\leq T}\|\rho-\rho_\infty\|_{L^2}^2\leq C\mathcal S_0^{\mathrm{nvac}}\rho_\infty^{1-\gamma},
\quad\sup_{0\leq t\leq T}\|p\|_{L^2}^2\leq C\mathcal S_0^{\mathrm{nvac}}\rho_\infty^{\gamma-1},\quad \|p_0\|_{L^2}^2
\leq C P_\infty\int G(\rho_0\mid\rho_\infty)dx.\label{eqeq3}
\end{gather}
\end{lemma}
\begin{proof}
It follows from the continuity equation \eqref{1}$_{1}$ that
\begin{align}\label{eq0}
\partial_tG(\rho\mid\rho_\infty)+\operatorname{div}\bigl(G(\rho\mid\rho_\infty)u\bigr)+p\operatorname{div} u=0.
\end{align}
Multiplying \eqref{1}$_{2}$ by $u$ and using \eqref{eq0}, we obtain
\begin{align}\label{eq1}
\frac{d}{d t}E(t)+\mu\|\nabla u\|_{L^2}^2+(\lambda+\mu)\|\operatorname{div} u\|_{L^2}^2=0.
\end{align}
Integrating \eqref{eq1}	over $(0,T)$ yields \eqref{eqba}.

For $\rho>0$, a direct calculation gives
\begin{align*}
\partial_\rho G(\rho\mid\rho_\infty)=\frac{A\gamma}{\gamma-1}\left(\rho^{\gamma-1}-\rho_\infty^{\gamma-1}\right),
\quad\partial_\rho^2G(\rho\mid\rho_\infty)=A\gamma\rho^{\gamma-2}.
\end{align*}
Since $G(\rho_\infty\mid\rho_\infty)=0$ and $\partial_\rho G(\rho_\infty\mid\rho_\infty)=0$,
Taylor's formula with integral remainder yields
\begin{align}\label{eqpo}
G(\rho\mid\rho_\infty)=A\gamma|\rho-\rho_\infty|^2\int_0^1(1-\theta)
\bigl((1-\theta)\rho_\infty+\theta\rho\bigr)^{\gamma-2} d\theta.
\end{align}
This identity also holds at $\rho=0$. By the upper bound of density and $R\geq1$,
\begin{align}\label{eqds}
(1-\theta)\rho_\infty\leq (1-\theta)\rho_\infty+\theta\rho\leq R\rho_\infty.
\end{align}
We distinguish the sign of $\gamma-2$. If $\gamma\geq2$, the function $s\mapsto s^{\gamma-2}$ is nondecreasing on $(0,\infty)$. Consequently,
\begin{align*}
\frac{\rho_\infty^{\gamma-2}}{\gamma}=\rho_\infty^{\gamma-2}\int_0^1(1-\theta)^{\gamma-1}\,  d\theta &\leq\int_0^1(1-\theta)\bigl((1-\theta)\rho_\infty+\theta\rho\bigr)^{\gamma-2} d\theta \\
&\leq( R\rho_\infty)^{\gamma-2}\int_0^1(1-\theta) d\theta =\frac{ (R\rho_\infty)^{\gamma-2}}{2},
\end{align*}
which along with \eqref{eqpo} indicates
\begin{align}\label{eq22}
A\rho_\infty^{\gamma-2}(\rho-\rho_\infty)^2\leq G(\rho\mid\rho_\infty)\leq
\frac{A\gamma}{2} (R\rho_\infty)^{\gamma-2}(\rho-\rho_\infty)^2.
\end{align}
If $1<\gamma<2$, the function $s\mapsto s^{\gamma-2}$ is decreasing on $(0,\infty)$. By \eqref{eqds}, we have
\begin{align*}
\frac{ (R\rho_\infty)^{\gamma-2}}{2}=(R\rho_\infty)^{\gamma-2}\int_0^1(1-\theta)\,  d\theta&\leq
\int_0^1(1-\theta)\bigl((1-\theta)\rho_\infty+\theta\rho\bigr)^{\gamma-2} d\theta\\
&\leq\rho_\infty^{\gamma-2}\int_0^1(1-\theta)^{\gamma-1}\,  d\theta=\frac{\rho_\infty^{\gamma-2}}{\gamma}.
\end{align*}
This together with \eqref{eqpo} implies
\begin{align*}
\frac{A\gamma}{2}(R\rho_\infty)^{\gamma-2}(\rho-\rho_\infty)^2\leq G(\rho\mid\rho_\infty)
\leq A\rho_\infty^{\gamma-2}(\rho-\rho_\infty)^2.
\end{align*}
Hence, \eqref{eqeq1} holds.
	
Finally, one infers from the fundamental theorem of calculus that
\begin{align*}
p=A\gamma(\rho-\rho_\infty)\int_0^1\bigl((1-\theta)\rho_\infty+\theta\rho\bigr)^{\gamma-1} d\theta.
\end{align*}
Since $\gamma-1>0$, it follows from \eqref{eqds} that
\begin{align*}
|p|\leq A\gamma (R\rho_\infty)^{\gamma-1}|\rho-\rho_\infty|,
\end{align*}
which combined with \eqref{eqeq1} indicates
\begin{align*}
p^2\leq A^2\gamma^2 (R\rho_\infty)^{2\gamma-2}(\rho-\rho_\infty)^2=A^2\gamma^2 R^{2\gamma-2}
\rho_\infty^\gamma\left(\rho_\infty^{\gamma-2}(\rho-\rho_\infty)^2\right)\leq C\rho_\infty^\gamma G(\rho\mid\rho_\infty),
\end{align*}
proving \eqref{eqeq2}. The desired \eqref{eqeq3} follows from \eqref{eqeq1} and \eqref{eqeq2}.
\end{proof}

\begin{lemma}\label{lembmo}
Let $0\leq\rho(x,t)\leq \bar\rho$. Then, for any $2\leq r<\infty$,
\begin{gather}
\|\nabla F\|_{L^r}+\|\nabla\omega\|_{L^r}\leq C\|\rho\dot u\|_{L^r},\label{eqer}\\
\|\nabla u\|_{L^4}^4\leq C\|\nabla u\|_{L^2}^2\left(\rho_\infty\|\sqrt\rho \dot u\|_{L^2}^2
+\rho_\infty^{2\gamma}\right),\label{L4}\\		
\|\nabla u\|_{L^3}^3\leq C\|\nabla u\|_{L^2}^2\big(\rho_\infty^{1/2}\|\sqrt\rho \dot u\|_{L^2}+\rho_\infty^\gamma\big).\label{L3}
\end{gather}
\end{lemma}
\begin{proof}
Taking the divergence and the curl of the momentum equation \eqref{1}$_{2}$ yields
\begin{align*}
\Delta F=\operatorname{div}(\rho\dot u),\quad\mu\Delta\omega=\operatorname{curl}(\rho\dot u).
\end{align*}
Then the standard elliptic estimates to the above systems give \eqref{eqer}.
	
We use a representation of the velocity gradient in terms of $F$, $\omega$, and $p$. Let
\begin{align*}
\mathcal R_j\triangleq\partial_j(-\Delta)^{-1/2}, \quad j=1,2,
\end{align*}
so that
\begin{align*}
\partial_i\partial_j(-\Delta)^{-1}=\mathcal R_i\mathcal R_j.
\end{align*}
Since $\omega=\partial_1u^2-\partial_2u^1$, we have
\begin{align*}
\Delta u^1=\partial_1\operatorname{div}u-\partial_2\omega,
\quad\Delta u^2=\partial_2\operatorname{div}u+\partial_1\omega.
\end{align*}
Combining these identities with
\begin{align*}
\operatorname{div}u=\frac{F+p}{2\mu+\lambda},
\end{align*}
we obtain that, for $i=1,2$,
\begin{align}\label{d2}
\partial_i u^1 =-\frac{1}{2\mu+\lambda}\mathcal R_i\mathcal R_1(F+p)+\mathcal R_i\mathcal R_2\omega,
\quad\partial_i u^2 =-\frac{1}{2\mu+\lambda}\mathcal R_i\mathcal R_2(F+p)-\mathcal R_i\mathcal R_1\omega.
\end{align}

We next estimate the BMO seminorm of $\nabla u$. Recall that
\begin{align*}
[f]_{\mathrm{BMO}}\triangleq\sup_Q\frac{1}{|Q|}\int_Q|f-f_Q| dx,
\quad f_Q\triangleq\frac{1}{|Q|}\int_Q f\,  dx,
\end{align*}	
where the supremum is taken over all squares $Q\subset\mathbb R^2$.	By the classical Calder\'o-Zygmund theory, the double Riesz transforms are bounded on  $\operatorname{BMO}\left(\mathbb{R}^{2}\right)$. Hence,
\begin{align*}
\left[\mathcal{R}_{i} \mathcal{R}_{j} f\right]_{\mathrm{BMO}} \leq C[f]_{\mathrm{BMO}},
\end{align*}
which together with \eqref{d2} gives
\begin{align*}
[\nabla u]_{\mathrm{BMO}}\leq C\left([F]_{\mathrm{BMO}}+[\omega]_{\mathrm{BMO}}
+[p]_{\mathrm{BMO}}\right).
\end{align*}
Using the Sobolev embedding $\dot H^1(\mathbb R^2)\hookrightarrow {\mathrm{BMO}}(\mathbb R^2)$
and the elementary estimate $[p]_{\mathrm{BMO}}\leq2\|p\|_{L^\infty}$, we deduce from \eqref{eqer} with $r=2$ that
\begin{align}\label{eqmo-density}
		[\nabla u]_{\mathrm{BMO}}
		 &\leq C\left(
		\|\nabla F\|_{L^2}
		+\|\nabla\omega\|_{L^2}
		+\|p\|_{L^\infty}
		\right)\notag\\
		& \leq C\left(
		\|\rho\dot u\|_{L^2}
		+\|p\|_{L^\infty}
		\right)\notag\\
& \leq C\big(R^{1/2}\rho_\infty^{1/2}\|\sqrt{\rho} \dot u\|_{L^2}+A\bigl( R^\gamma+1\bigr)\rho_\infty^\gamma\big)\notag\\
& \leq C\big(\rho_\infty^{1/2}\|\sqrt{\rho} \dot u\|_{L^2}+\rho_\infty^\gamma\big).
\end{align}

By Lemma \ref{ggn} and the embedding $L^2(\mathbb R^2)\hookrightarrow L^{2,\infty}(\mathbb R^2)$, one gets that, for every $f\in L^2(\mathbb R^2)\cap {\mathrm{BMO}}(\mathbb R^2)$,
\begin{align}\label{eq:BMO-interpolation}
\|f\|_{L^r}\leq C \|f\|_{L^2}^{2/r}[f]_{\mathrm{BMO}}^{1-2/r},\quad 2<r<\infty.
\end{align}
Applying this inequality componentwise to $\nabla u$ with $r=4$ and using \eqref{eqmo-density}, we obtain
	\begin{align*}
		\|\nabla u\|_{L^4}^4
		 \leq C\|\nabla u\|_{L^2}^2
		[\nabla u]_{\mathrm{BMO}}^2 \leq C\|\nabla u\|_{L^2}^2
		\left(
		\rho_\infty
		\|\sqrt{\rho} \dot u\|_{L^2}^2
		+\rho_\infty^{2\gamma}
		\right),
	\end{align*}
proving \eqref{L4}.
	Similarly, taking $r=3$ in \eqref{eq:BMO-interpolation} gives \eqref{L3}.
\end{proof}

\begin{remark}\label{u1}
From \eqref{d2} and the properties of the Riesz transform, one has that
\begin{align}\label{eqho}
\|\nabla u\|_{L^r}\leq C\left(\|F\|_{L^r}+\|\omega\|_{L^r}+\|p\|_{L^r}\right)
\end{align}
for any $1<r<\infty$. Unlike \eqref{eqho}, the estimate \eqref{L4} retains the factor $\|\nabla u\|_{L^2}^2$ in the pressure contribution. This structure is particularly useful because the basic energy inequality provides control of $\int_0^T\|\nabla u\|_{L^2}^2dt$, whereas it does not directly yield a uniform-in-time bound of $\int_0^T\|p\|_{L^4}^4dt$.
\end{remark}

\subsection{\textit{A priori} estimates I: scaling invariant estimates}
	
\begin{proposition}[Bootstrap argument]\label{boot}
Let $\mathcal S_{t}^{\mathrm{nvac}}\triangleq\rho_\infty^{1-\gamma}\left(\|\nabla u(t)\|_{L^2}^2+P_\infty E(t)\right)$. Assume that
\begin{align}\label{cp1}
0\leq\rho(x,t)\leq2\bar\rho,\quad\mathcal S_{t}^{\mathrm{nvac}}\leq2\varepsilon
\end{align}
for all $(x,t)\in \mathbb R^3\times[0,T]$. Then there exist positive constants
$\varepsilon_1=\varepsilon_1(\RR,A,\gamma,\mu,\lambda)$ and $c_0=c_0(\RR,A,\gamma,\mu,\lambda)$
such that
\begin{align*}
0\leq\rho(x,t)\leq\frac32\bar\rho,\quad\mathcal S_{t}^{\mathrm{nvac}}\leq\frac32\varepsilon
\end{align*}
provided
\begin{align}\label{cp2}
c_0\mathcal S^{\mathrm{nvac}}_0\leq\varepsilon\leq\varepsilon_1.
\end{align}
\end{proposition}

\begin{lemma}\label{3lem1}
Under the assumption \eqref{cp1}, it holds that
\begin{align}\label{eqlow}
\sup_{0\leq t\leq T}\mathcal S_{t}^{\mathrm{nvac}}+\rho_\infty^{1-\gamma}
\int_0^T\|\sqrt\rho \dot u\|_{L^2}^2\, dt\leq C\mathcal S_0^{\mathrm{nvac}}.
\end{align}
\end{lemma}
\begin{proof}
Multiplying \eqref{1}$_{2}$ by $\dot u$ and integrating the resultant over $\mathbb{R}^2$, we derive from \eqref{L3} and Young's inequality that
\begin{align}\label{eqf1}
&\frac{d}{dt}\left[\frac\mu2\|\nabla u\|_{L^2}^2+\frac{\mu+\lambda}{2}\|\operatorname{div} u\|_{L^2}^2
-\int p\operatorname{div} u\, dx\right]+\|\sqrt\rho \dot u\|_{L^2}^2\notag\\
&=-\mu\int\partial_i u^j\partial_i u^k\partial_k u^j dx+\frac\mu2\int \operatorname{div} u\,|\nabla u|^2 dx
-(\mu+\lambda)\int\operatorname{div} u\,\partial_j u^k\partial_k u^j dx\notag\\
&\quad+\frac{\mu+\lambda}{2}\int (\operatorname{div} u)^3 dx+\int p\,\partial_j u^k\partial_k u^j dx+(\gamma-1)\int p(\operatorname{div} u)^2 dx+\gamma P_\infty\|\operatorname{div} u\|_{L^2}^2\notag\\
&\leq C\|\nabla u\|_{L^3}^3+C\bar{\rho}^\gamma\|\nabla u\|_{L^2}^2+C\rho_\infty^\gamma\|\nabla u\|_{L^2}^2\notag \\
 	&\leq C\rho_\infty^{1/2}\|\nabla u\|_{L^2}^2
 	\|\sqrt\rho \dot u\|_{L^2}
 	+C\rho_\infty^\gamma\|\nabla u\|_{L^2}^2\notag \\
 	&\leq\frac12\|\sqrt\rho \dot u\|_{L^2}^2
 	+C\rho_\infty\|\nabla u\|_{L^2}^4
 	+C\rho_\infty^\gamma\|\nabla u\|_{L^2}^2.
 \end{align}
Hence
\begin{align}\label{eqf11}
\frac{d}{dt}\left[
	\frac{\mu}{2}\|\nabla u\|_{L^2}^2
	+\frac{\mu+\lambda}{2}\|\operatorname{div} u\|_{L^2}^2
	-\int p\operatorname{div} u\, dx
	\right]
	+\frac12\|\sqrt\rho \dot u\|_{L^2}^2
 \leq C\rho_\infty\|\nabla u\|_{L^2}^4
	+C\rho_\infty^\gamma\|\nabla u\|_{L^2}^2.
\end{align}
	
Introduce the modified energy functional
\begin{align}
		\mathcal H=\mathcal H(t)
		\triangleq\frac\mu2\|\nabla u(t)\|_{L^2}^2
		+\frac{\mu+\lambda}{2}\|\operatorname{div} u(t)\|_{L^2}^2
		-\int p(t)\operatorname{div} u(t)\, dx+K\rho_\infty^\gamma E(t)
		\label{eqmo}
	\end{align}
with a positive constant $K=K(\RR,A,\gamma,\mu,\lambda)$.
If $K$ is sufficiently large, one gets from Young's inequality and \eqref{eqeq2} that
	\begin{align}\label{eqde}
		c\left[\|\nabla u(t)\|_{L^2}^2
		+\rho_\infty^\gamma E(t)\right]
		\leq\mathcal H
		\leq C\left[\|\nabla u(t)\|_{L^2}^2
		+\rho_\infty^\gamma E(t)\right].	
	\end{align}
Adding \eqref{eq1} multiplied by $K\rho_\infty^\gamma$ to \eqref{eqf11} and increasing $K$ if necessary yields
\begin{align}\label{zui}
		\frac{d}{dt}\mathcal H
		+\|\sqrt\rho \dot u\|_{L^2}^2
		+\rho_\infty^\gamma\|\nabla u\|_{L^2}^2
		\leq C\rho_\infty\|\nabla u\|_{L^2}^4\leq C\rho_\infty\|\nabla u\|_{L^2}^2\mathcal H.
	\end{align}
It follows from \eqref{eqba} that
\begin{align*}
\rho_\infty\int_0^T\|\nabla u\|_{L^2}^2\, dt\leq C\mathcal S_0^{\mathrm{nvac}}\leq C \varepsilon,
\end{align*}
which together with \eqref{zui}, the fact $\mathcal H(0)\leq C\mathcal S_0^{ \mathrm{nvac}}\rho_\infty^{\gamma-1}$, and Gronwall's inequality leads to
\begin{align*}
\sup_{0\leq t\leq T}\mathcal H(t)+\int_0^T\|\sqrt\rho \dot u\|_{L^2}^2 dt
\leq C\mathcal S_0^{\mathrm{nvac}}\rho_\infty^{\gamma-1}.
\end{align*}
This along with \eqref{eqde} shows \eqref{eqlow}.
\end{proof}

\begin{lemma}\label{lem:material}
Under the assumption \eqref{cp1}, it holds that
\begin{align}\label{w0}
\sup_{0\leq t\leq T}\left( \min\{1,\rho_\infty^\gamma t\}\|\sqrt\rho \dot u\|_{L^2}^2\right)
+\int_0^T\min\{1,\rho_\infty^\gamma t\}\|\nabla\dot u\|_{L^2}^2dt
\leq C\mathcal S_0^{\mathrm {nvac}}\rho_\infty^{2\gamma-1}.
\end{align}
\end{lemma}
\begin{proof}
Applying $\partial_t+\partial_k(u^k\,\cdot)$ to the $j$th component of the momentum equation \eqref{1}$_{2}$, we obtain\footnote{See Appendix A for the derivation.}
\begin{align}\label{c19}
&\rho(\partial_t\dot u^j+u\cdot\nabla\dot u^j)-\mu\Delta\dot u^j-(\mu+\lambda)\partial_j\operatorname{div}\dot u\notag\\
&=\mu\partial_i\big[(\operatorname{div} u)\partial_i u^j-\partial_i u^k\partial_k u^j\big]
-\mu\partial_k\big(\partial_i u^k\partial_i u^j\big)+(\mu+\lambda)\partial_j
\big[(\operatorname{div} u)^2-\partial_i u^k\partial_k u^i\big]\notag\\
&\quad-(\mu+\lambda)\partial_k\big[(\operatorname{div} u)\partial_j u^k\big]
+(\gamma-1)\partial_j(P\operatorname{div} u)
+\partial_k(P\partial_j u^k).
\end{align}
Multiplying \eqref{c19} by $\dot u^j$, summing over $j$, and integrating by parts over $\mathbb R^2$, we deduce from H\"older's inequality, Young's inequality, and \eqref{L4} that
\begin{align*}
&\frac12\frac{d}{dt}\|\sqrt\rho \dot u\|_{L^2}^2+\mu\|\nabla\dot u\|_{L^2}^2+(\mu+\lambda)\|\operatorname{div}\dot u\|_{L^2}^2\\
&=\mu\int\big[\partial_i u^k\partial_k u^j-(\operatorname{div} u)\partial_i u^j\big]
\partial_i\dot u^j dx+\mu\int\partial_i u^k\partial_i u^j\partial_k\dot u^j dx\\
&\quad+(\mu+\lambda)\int\big[\partial_i u^k\partial_k u^i-(\operatorname{div} u)^2\big]
\operatorname{div}\dot u dx+(\mu+\lambda)\int(\operatorname{div} u)\partial_j u^k\partial_k\dot u^j dx\\
&\quad-(\gamma-1)\int P\operatorname{div} u\operatorname{div}\dot u dx
-\int P\partial_j u^k\partial_k\dot u^j dx\\
&\leq C\|\nabla u\|_{L^4}^2\|\nabla\dot u\|_{L^2}
+C\rho_\infty^\gamma\|\nabla u\|_{L^2}\|\nabla\dot u\|_{L^2}\\
&\leq \frac\mu2\|\nabla\dot u\|_{L^2}^2+C\rho_\infty\|\nabla u\|_{L^2}^2\|\sqrt\rho \dot u\|_{L^2}^2+C\rho_\infty^{2\gamma}\|\nabla u\|_{L^2}^2.
\end{align*}
Thus,
\begin{align}\label{w1}
\frac{d}{dt}\|\sqrt\rho\dot u\|_{L^2}^2
+\|\nabla\dot u\|_{L^2}^2\leq C\rho_\infty\|\nabla u\|_{L^2}^2\|\sqrt\rho\dot u\|_{L^2}^2
+C\rho_\infty^{2\gamma}\|\nabla u\|_{L^2}^2.
\end{align}
Multiplying \eqref{w1} by $\min\{1,\rho_\infty^\gamma t\}$, integrating the resultant in $[0,T]$, and using \eqref{eqlow}, one has that
\begin{align*}
&\sup_{0\leq t\leq T}\min\{1,\rho_\infty^\gamma t\}\|\sqrt\rho \dot u\|_{L^2}^2
+\int_0^T\min\{1,\rho_\infty^\gamma t\}\|\nabla\dot u\|_{L^2}^2 dt\\
&\leq C\rho_\infty^\gamma\int_0^{\min\{\rho_\infty^{-\gamma},T\}}
\|\sqrt\rho \dot u\|_{L^2}^2dt+C\rho_\infty\sup_{0\leq t\leq T}\|\nabla u\|_{L^2}^2
\int_0^T\|\sqrt\rho \dot u\|_{L^2}^2 dt+C\rho_\infty^{2\gamma}\int_0^T\|\nabla u\|_{L^2}^2 dt\\
&\leq C\mathcal S_0^{\mathrm{nvac}}\rho_\infty^{2\gamma-1}+C(\mathcal S_0^{ \mathrm{nvac}})^2\rho_\infty^{2\gamma-1}\\
&\leq C\mathcal S_0^{\mathrm{nvac}}\rho_\infty^{2\gamma-1},
\end{align*}
proving \eqref{w0}.
\end{proof}

\begin{remark}[Choice of the time weight]\label{weight}
The weight $\min\{1,\rho_\infty^\gamma t\}$ in \eqref{w0} is chosen so that the contribution from its time derivative has the same power of $\rho_\infty$ as the other terms in the weighted estimate. This yields a bound whose constant is independent of $\rho_\infty$.
	
To explain this choice, let $\tau>0$ be an arbitrary cutoff time and set
\begin{align*}
\sigma_\tau(t)=\min\left\{1,\frac{t}{\tau}\right\}.
\end{align*}
Multiplying \eqref{w1} by $\sigma_\tau$ and integrating in time gives
\begin{align*}
&\sup_{0\leq t\leq T}\sigma_\tau(t)\|\sqrt\rho \dot u(t)\|_{L^2}^2+\int_0^T\sigma_\tau(t)\|\nabla\dot u(t)\|_{L^2}^2dt\notag\\
&\leq\frac{C}{\tau}\int_0^{\min\{\tau,T\}}\|\sqrt\rho \dot u\|_{L^2}^2\,dt+C\rho_\infty\int_0^T\sigma_\tau(t)\|\nabla u\|_{L^2}^2\|\sqrt\rho \dot u\|_{L^2}^2dt
+C\rho_\infty^{2\gamma}\int_0^T\sigma_\tau(t)\|\nabla u\|_{L^2}^2dt.
\end{align*}
	
Recall the lower-order estimates
\begin{align*}
\sup_{0\leq t\leq T}\|\nabla u(t)\|_{L^2}^2\leq C\varepsilon\rho_\infty^{\gamma-1},
\quad\int_0^T\|\sqrt\rho \dot u\|_{L^2}^2dt\leq C\varepsilon\rho_\infty^{\gamma-1},
\quad\int_0^T\|\nabla u\|_{L^2}^2dt\leq C\varepsilon\rho_\infty^{-1}.	
\end{align*}
The term arising from the derivative of the weight therefore satisfies
\begin{align*}
\frac{1}{\tau}\int_0^{\min\{\tau,T\}}\|\sqrt\rho \dot u\|_{L^2}^2dt
\leq C\varepsilon\tau^{-1}\rho_\infty^{\gamma-1}.
\end{align*}
Since $\sigma_\tau\leq1$, the remaining two terms are estimated as follows:
\begin{gather*}
C\rho_\infty\int_0^T\sigma_\tau(t)\|\nabla u\|_{L^2}^2\|\sqrt\rho \dot u\|_{L^2}^2\,dt
\leq C\rho_\infty\sup_{0\leq t\leq T}\|\nabla u(t)\|_{L^2}^2
\int_0^T\|\sqrt\rho \dot u\|_{L^2}^2dt\leq C\varepsilon^2\rho_\infty^{2\gamma-1},
\\ C\rho_\infty^{2\gamma}\int_0^T\sigma_\tau(t)\|\nabla u\|_{L^2}^2dt
\leq C\varepsilon\rho_\infty^{2\gamma-1}.
\end{gather*}
Consequently, for $0<\varepsilon\leq1$,
\begin{align}\label{GG}
\sup_{0\leq t\leq T}\min\left\{1,\frac{t}{\tau}\right\}\|\sqrt\rho \dot u(t)\|_{L^2}^2
+\int_0^T\min\left\{1,\frac{t}{\tau}\right\}\|\nabla\dot u\|_{L^2}^2\,dt
\leq C\varepsilon\left(\tau^{-1}\rho_\infty^{\gamma-1}+\rho_\infty^{2\gamma-1}\right).
\end{align}
Now we choose $\tau$ so that the two terms on the right-hand side have the same dependence on $\rho_\infty$. This gives
\begin{align*}
\tau=\rho_\infty^{-\gamma},\quad\sigma_\tau(t)=\min\{1,\rho_\infty^\gamma t\}.
\end{align*}
With this choice, \eqref{GG} becomes
\begin{align}\label{JJ}
\sup_{0\leq t\leq T}\min\{1,\rho_\infty^\gamma t\}\|\sqrt\rho \dot u(t)\|_{L^2}^2+\int_0^T\min\{1,\rho_\infty^\gamma t\}\|\nabla\dot u\|_{L^2}^2\,dt\leq C\varepsilon\rho_\infty^{2\gamma-1}.
\end{align}
	
If the fixed cutoff $\min\{1,t\}$ is used instead, then $\tau=1$, and \eqref{GG} yields
\begin{align*}
\sup_{0\leq t\leq T}\min\{1,t\}\|\sqrt\rho \dot u(t)\|_{L^2}^2+\int_0^T\min\{1,t\}\|\nabla\dot u\|_{L^2}^2\,dt
\leq C\varepsilon\left(\rho_\infty^{\gamma-1}+\rho_\infty^{2\gamma-1}\right)
=C\varepsilon\rho_\infty^{2\gamma-1}\left(1+\rho_\infty^{-\gamma}\right).
\end{align*}
For arbitrary $\rho_\infty>0$, the factor $1+\rho_\infty^{-\gamma}$ cannot be absorbed into a constant
independent of $\rho_\infty$.
	
The chosen weight also serves the purpose  in the subsequent analysis. The weight equals $1$ for $t\geq\rho_\infty^{-\gamma}$, so that the weighted estimate directly provides an unweighted integral bound on the
long-time interval, without requiring an additional decay estimate.
	
More precisely, \eqref{JJ} implies
\begin{align*}
\int_0^{\min\{\rho_\infty^{-\gamma},T\}}t\|\nabla\dot u\|_{L^2}^2dt\leq C\varepsilon\rho_\infty^{\gamma-1},
\end{align*}
and, whenever $T>\rho_\infty^{-\gamma}$,
\begin{align*}
\int_{\rho_\infty^{-\gamma}}^T\|\nabla\dot u\|_{L^2}^2dt\leq C\varepsilon\rho_\infty^{2\gamma-1}.
\end{align*}
These powers of $\rho_\infty$ match the factors arising from spatial interpolation and H\"older's inequality in time in the subsequent effective-viscous-flux estimates. Together with the other a priori bounds, they yield a small upper bound, independent of $\rho_\infty$, for
\begin{align*}
\int_0^{\min\{\rho_\infty^{-\gamma},T\}}\|F\|_{L^\infty}\,dt.
\end{align*}
This uniformity is needed to obtain a smallness threshold independent of $\rho_\infty$ when applying Zlotnik's inequality to control the density.
\end{remark}

\subsection{Proof of Proposition \ref{boot}}\label{sub}

\begin{proof}[Proof of Proposition \ref{boot}]
{\it Step 1: Closing the bootstrap estimate for $\mathcal S_t^{\mathrm{nvac}}$.}
By \eqref{eqlow}, we choose
\begin{align}\label{c0}
c_0\triangleq\max\left\{\frac12,\frac{2C}{3}\right\}.
\end{align}
Since $c_0\mathcal S_{0}^{\mathrm{nvac}}\leq\varepsilon$, it follows that
\begin{align*}
\mathcal S_t^{\mathrm{nvac}}\leq C \mathcal S_0^{\mathrm{nvac}}
\leq\frac{C}{c_0}\varepsilon\leq\frac32\varepsilon.
\end{align*}

{\it Step 2: Closing the bootstrap estimate for density.}
We infer from Gagliardo--Nirenberg inequality, \eqref{eqd2}, \eqref{eqeq3}, \eqref{eqer}, \eqref{cp1}, and \eqref{eqlow} that
\begin{align}\label{x1}
\|F\|_{L^\infty}&\leq C\|F\|_{L^2}^{1/3}\|\nabla F\|_{L^4}^{2/3}\notag\\
&\leq [C\varepsilon^{1/2}\rho_\infty^{(\gamma-1)/2}]^{1/3}[\|\rho\dot u\|_{L^4}]^{2/3}\notag\\
&\leq [C\varepsilon^{1/2}\rho_\infty^{(\gamma-1)/2}]^{1/3}\big[C\rho_\infty^{3/4}
\|\sqrt\rho \dot u\|_{L^2}^{1/2}\|\nabla\dot u\|_{L^2}^{1/2}
+C\rho_\infty^{1/2}\|\rho-\rho_\infty\|_{L^2}^{1/2}\|\nabla\dot u\|_{L^2}\big]^{2/3}\notag\\
&\leq C\varepsilon^{1/6}\rho_\infty^{(\gamma+2)/6}
\|\sqrt\rho \dot u\|_{L^2}^{1/3}\|\nabla\dot u\|_{L^2}^{1/3}
+C\varepsilon^{1/3}\rho_\infty^{1/3}\|\nabla\dot u\|_{L^2}^{2/3},
\end{align}
where we have used
\begin{gather*}
\|F\|_{L^2}\leq C\left(\|\nabla u\|_{L^2}+\|p\|_{L^2}\right)\leq C\varepsilon^{1/2}\rho_\infty^{(\gamma-1)/2},\\
\|\nabla F\|_{L^4}\leq C\|\rho\dot u\|_{L^4}
\leq C\rho_\infty\|\dot u\|_{L^4}\leq C\rho_\infty^{3/4}
\|\sqrt\rho \dot u\|_{L^2}^{1/2}\|\nabla\dot u\|_{L^2}^{1/2}
+C\rho_\infty^{1/2}\|\rho-\rho_\infty\|_{L^2}^{1/2}\|\nabla\dot u\|_{L^2}.
\end{gather*}
	
If $0\leq t\leq\min\{\rho_\infty^{-\gamma},T\}$, the time weight in \eqref{w0} equals $\rho_\infty^\gamma t$.
Consequently,
\begin{align*}
\int_0^{\min\{\rho_\infty^{-\gamma},T\}}t\|\nabla\dot u\|_{L^2}^2\, dt
\leq C\varepsilon\rho_\infty^{\gamma-1},
\end{align*}
which together with H\"older's inequality and \eqref{eqlow} yields
\begin{align*}
&\int_0^{\min\{\rho_\infty^{-\gamma},T\}}\|\sqrt\rho \dot u\|_{L^2}^{1/3}
\|\nabla\dot u\|_{L^2}^{1/3}\, dt\\
&\leq\left(\int_0^{\min\{\rho_\infty^{-\gamma},T\}}
\|\sqrt\rho \dot u\|_{L^2}^2\, dt\right)^{1/6}
\left(\int_0^{\min\{\rho_\infty^{-\gamma},T\}}
t\|\nabla\dot u\|_{L^2}^2\, dt\right)^{1/6}
\left(\int_0^{\min\{\rho_\infty^{-\gamma},T\}}t^{-1/4}\, dt\right)^{2/3}\\[2mm]
&\leq C\left(\varepsilon\rho_\infty^{\gamma-1}\right)^{1/3}\rho_\infty^{-\gamma/2}
=C\varepsilon^{1/3}\rho_\infty^{-(\gamma+2)/6}.
\end{align*}
Similarly,
\begin{align*}
\int_0^{\min\{\rho_\infty^{-\gamma},T\}}\|\nabla\dot u\|_{L^2}^{2/3}\, dt
&\leq\left(\int_0^{\min\{\rho_\infty^{-\gamma},T\}}t\|\nabla\dot u\|_{L^2}^2\, dt\right)^{1/3}
\left(\int_0^{\min\{\rho_\infty^{-\gamma},T\}}t^{-1/2}\, dt\right)^{2/3}\\[2mm]
&\leq C\left(\varepsilon\rho_\infty^{\gamma-1}\right)^{1/3}\rho_\infty^{-\gamma/3}
=C\varepsilon^{1/3}\rho_\infty^{-1/3}.
\end{align*}
Integrating \eqref{x1} and using the above two estimates shows
\begin{align}\label{eqsh}
\int_0^{\min\{\rho_\infty^{-\gamma},T\}}\|F\|_{L^\infty} dt
\leq C\big(\varepsilon^{1/2}+\varepsilon^{2/3}\big).
\end{align}
If $T>\rho_\infty^{-\gamma}$, \eqref{w0} gives
\begin{align*}
\int_{\rho_\infty^{-\gamma}}^T\|\nabla\dot u\|_{L^2}^2dt
\leq C\varepsilon\rho_\infty^{2\gamma-1}.
\end{align*}
This combined with \eqref{x1}, \eqref{eqlow}, and Cauchy--Schwarz inequality indicates
\begin{align}\label{eqsh1}
\int_{\rho_\infty^{-\gamma}}^T\|F\|_{L^\infty}^3 dt
&\leq C\varepsilon^{1/2}\rho_\infty^{(\gamma+2)/2}
\left(\int_{\rho_\infty^{-\gamma}}^T\|\sqrt\rho \dot u\|_{L^2}^2 dt\right)^{1/2}
\left(\int_{\rho_\infty^{-\gamma}}^T\|\nabla\dot u\|_{L^2}^2 dt\right)^{1/2}+C\varepsilon\rho_\infty
\int_{\rho_\infty^{-\gamma}}^T\|\nabla\dot u\|_{L^2}^2 dt\notag\\
&\leq C\rho_\infty^{2\gamma}\big(\varepsilon^{3/2}+\varepsilon^2\big).
\end{align}

We next apply Lemma \ref{zo} to the density. For each $x_0\in\mathbb R^2$, let $x(t)$ solve
\begin{align*}
\frac{d}{dt}x(t)=u(x(t),t),\quad x(0)=x_0.
\end{align*}
Set
\begin{align*}	
y(t)=\rho(x(t),t), \quad b(t) =-\frac{1}{2\mu+\lambda}
\int_0^t\rho(x(s),s)F(x(s),s)\,  ds.
\end{align*}
Since $F=(2\mu+\lambda)\operatorname{div} u-A(\rho^\gamma-\rho_\infty^\gamma)$, one gets from the continuity equation \eqref{1}$_{1}$ that
\begin{align*}
y'(t)=g(y(t))+b'(t),\quad g(s)=-\frac{As}{2\mu+\lambda}
\left(s^\gamma-\rho_\infty^\gamma\right),
\quad s\geq0.
\end{align*}
Choose
\begin{align*}
\zeta=\frac{5\bar\rho}{4},\quad\delta=\frac A2\left[
\left(\frac{5\bar\rho}{4}\right)^\gamma-\rho_\infty^\gamma\right]>0,
\quad N_1=\frac{2\bar\rho\,\delta}{2\mu+\lambda}.
\end{align*}
Noting that for $a\geq0$,
\begin{align*}
a=\delta^{2/3}(\delta^{-2/3}a)
\leq\frac{2\delta}{3}+\frac{(\delta^{-2/3}a)^3}{3}
\leq\delta+\delta^{-2}a^3,
\end{align*}
and hence, for any $0\leq t_1<t_2\leq T$,
	\begin{align*}
		\int_{t_1}^{t_2}\|F\|_{L^\infty}\, dt
		&\leq
		\int_{(t_1,t_2)\cap(0,\rho_\infty^{-\gamma})}
		\|F\|_{L^\infty}\, dt
		+\delta(t_2-t_1)+\delta^{-2}
		\int_{(t_1,t_2)\cap(\rho_\infty^{-\gamma},T)}
		\|F\|_{L^\infty}^3\, dt\\
		&\leq C\big(\varepsilon^{1/2}+\varepsilon^{2/3}\big)
		+\delta(t_2-t_1)+C\delta^{-2}\rho_\infty^{2\gamma}
		\big(\varepsilon^{3/2}+\varepsilon^2\big),
	\end{align*}
which together with \eqref{eqsh}, \eqref{eqsh1}, and the bootstrap bound $0\leq\rho\leq2\bar\rho$ concludes that
\begin{align*}
|b(t_2)-b(t_1)|\leq\frac{2\bar\rho}{2\mu+\lambda}
\int_{t_1}^{t_2}\|F\|_{L^\infty}\, dt\leq N_0+N_1(t_2-t_1),
\end{align*}
where
\begin{align}\label{eqN}
N_0=\frac{C\bar\rho}{2\mu+\lambda}
\big[\varepsilon^{1/2}+\varepsilon^{2/3}+\delta^{-2}\rho_\infty^{2\gamma}
\big(\varepsilon^{3/2}+\varepsilon^2\big)\big].
\end{align}
For every $s\geq\zeta$, we also have
	\begin{align*}
	g(s)\leq-\frac{A\zeta}{2\mu+\lambda}
	\left(\zeta^\gamma-\rho_\infty^\gamma\right)
	=-\frac{2\zeta\delta}{2\mu+\lambda}
	=-\frac54N_1\leq-N_1.
\end{align*}
Therefore, all the hypotheses of Zlotnik's inequality are satisfied.	
To choose the smallness threshold independently of $\rho_\infty$, observe that
	\begin{align*}
	\delta^{-2}\rho_\infty^{2\gamma}
	=\frac{4}{A^2\left[(5R/4)^\gamma-1\right]^2}.
\end{align*}
It follows from \eqref{eqN} that
\begin{align*}
\frac{N_0}{\bar\rho}=\frac{C}{2\mu+\lambda}\left[\varepsilon^{1/2}+\varepsilon^{2/3}
+\frac{4\left(\varepsilon^{3/2}+\varepsilon^2\right)}{A^2\left[(5R/4)^\gamma-1\right]^2}\right].
\end{align*}
Choose $\varepsilon_1$ sufficiently small such that $N_0\leq\bar\rho/4$.
Since $y(0)\leq\bar\rho<\zeta$, Zlotnik's inequality yields
\begin{align*}
\rho(x(t),t)\leq\frac{5\bar\rho}{4}+N_0\leq\frac32\bar\rho. \tag*{\qedhere}
\end{align*}
\end{proof}

\subsection{Proof of Theorem \ref{th1}}

Under the condition \eqref{d1}, a unique local strong solution $(\rho,u)$ to the problem \eqref{1}--\eqref{3} verifying \eqref{def1} has been obtained in \cite{Luo2012}. Let $T_{\max}>0$ denote its maximal existence time. By Lemma \ref{blowup}, if $T_{\max}<\infty$, then
\begin{align}\label{bb1}
\lim_{T\rightarrow T_{\max}}\|\rho\|_{L^\infty(0,T;L^\infty)}=\infty.
\end{align}
Furthermore, Proposition \ref{boot} together with a standard continuity argument shows that the density bound in \eqref{cp1} holds on $[0,T]$ for every $T<T_{\max}$ provided \eqref{S0} holds.
This contradicts \eqref{bb1}. Therefore, $T_{\max}=\infty$. \hfill$\square$

\section{Compressible Navier--Stokes equations with far-field vacuum}\label{sec4}

In this section we consider the Cauchy problem \eqref{1}--\eqref{3} with far-field vacuum. For initial data satisfying \eqref{ini}, the local well-posedness of strong solutions has been derived in \cite{LiLiang2014}. Throughout this section, we assume that $(\rho,u)$ is the strong solution to \eqref{1}--\eqref{3} on $\mathbb R^2\times[0,T]$ for some $T>0$. Our main objective is to derive the {\it a priori} estimates independent of $T$. Unless otherwise specified, $C$ denotes a generic positive constant depending only on $\mu$, $\lambda$, and $\gamma$, whose value may vary from line to line. In particular, $C$ is independent of $T$, $m_0$, $\bar\rho$, and the norms of initial data.

In what follows, we first derive basic energy and elliptic estimates and then establish a lower-order bootstrap estimate. After introducing the intrinsic scale and the first-moment parameters, we state the main bootstrap proposition. Time-weighted material-derivative estimates and mass localization provide the effective viscous flux bounds needed to close the density upper bound.

\subsection{Basic energy and auxiliary estimates}
Integrating \eqref{1}$_{1}$ gives
\begin{align}\label{eqma}
\int\rho(x,t)\, d x=m_0= \int \rho_0\,dx.
\end{align}
For later use, we recall
\begin{gather*}
D(t)=\|\nabla u(t)\|_{L^{2}}^{2}+\|P(t)\|_{L^{2}}^{2},\quad D_{0}=D(0),\\
E(t)=\frac12 \|\sqrt\rho u(t)\|_{L^{2}}^{2}+\frac1{\gamma-1}\|P(t)\|_{L^{1}},
\quad E_0= E(0),\\
 \mathcal S_0\triangleq\bar\rho\bigl(E_0+m_0D_{0}\bigr),\quad
F =(2\mu+\lambda)\operatorname{div} u-P.
\end{gather*}
Throughout this subsection, we assume that
\begin{align*}
0\leq\rho(x,t)\leq 2\bar\rho
\qquad\text{for all }(x,t)\in\R^2\times[0,T].
\end{align*}
This density bound
will later be part of the unified bootstrap hypothesis.
\begin{lemma}[Basic energy estimates]
It holds that
\begin{align}\label{eqhe}
\sup_{0\le t\le T}E(t)+\int_0^T\left(\mu\|\nabla u\|_{L^{2}}^{2}
+(\mu+\lambda)\|\operatorname{div} u\|_{L^{2}}^{2}\right)dt\leq CE_0.
\end{align}
\end{lemma}
\begin{proof}
Testing \eqref{1}$_{2}$ by $u$ and using \eqref{1}$_{1}$, we obtain	\eqref{eqhe}.
\end{proof}

\begin{lemma}
It holds that
	\begin{gather}
		\|F\|_{L^{3}}^3+\|\omega\|_{L^{3}}^3
		\le C\bar\rho^{1/2}\|\sqrt\rho \dot u\|_{L^{2}}D(t),
		\label{eq31}\\[1mm]
		\|F\|_{L^{4}}^4+\|\omega\|_{L^{4}}^4
		\le C\bar\rho D(t)\|\sqrt\rho \dot u\|_{L^{2}}^{2},
		\label{eqe33}\\[1mm]
		\|\nabla u\|_{L^{3}}^3
		\le C\bar\rho^{1/2}\|\sqrt\rho \dot u\|_{L^{2}}D(t)
		+C\|P\|_{L^{3}}^3,
		\label{eq32}\\[1mm]
		\|\nabla u\|_{L^{4}}^4
		\le C\bar\rho D(t)\|\sqrt\rho \dot u\|_{L^{2}}^{2}
		+C\|P\|_{L^{4}}^4.
		\label{eq34}
	\end{gather}
\end{lemma}

\begin{proof}
	It follows form \eqref{eqer} and \eqref{eqho} that
	\begin{gather}\label{L2}
		\|\nabla F\|_{L^{2}}+\|\nabla\omega\|_{L^{2}}
		\le C\bar\rho^{1/2}\|\sqrt\rho \dot u\|_{L^{2}},\\
		\|\nabla u\|_{L^{3}}
		\le C\bigl(\|F\|_{L^{3}}+\|\omega\|_{L^{3}}+\|P\|_{L^{3}}\bigr),
		\label{eq35}\\
		\|F\|_{L^{2}}+\|\omega\|_{L^{2}}\le C D(t)^{1/2}.\label{t2}
\end{gather}
The Gagliardo--Nirenberg inequality $\|f\|_{L^{3}}^3\le C\|f\|_{L^{2}}^{2}\|\nabla f\|_{L^{2}}$
applied to $F$ and $\omega$ gives \eqref{eq31}, while the desired \eqref{eq32} follows from \eqref{eq35} and \eqref{eq31}. Similarity, one deduces from
\begin{align*}
\|f\|_{L^{4}}^4\le C\|f\|_{L^{2}}^{2}\|\nabla f\|_{L^{2}}^{2}
\end{align*}
that \eqref{eqe33} and \eqref{eq34} holds.
\end{proof}

\begin{lemma}[Pressure estimate with mass dependence]
It holds that
\begin{gather}
\|P\|_{L^{2}}\le C m_0^{1/2}\|\sqrt\rho \dot u\|_{L^{2}}+C\|\nabla u\|_{L^{2}},
\label{eqp1}\\
\int_0^T D(t) d t\le C E_0+C m_0\int_0^T\|\sqrt\rho \dot u\|_{L^{2}}^{2} d t.\label{eqp2}
\end{gather}
\end{lemma}
\begin{proof}
Taking the divergence of \eqref{1}$_{2}$ yields
\begin{align}\label{eqp3}
P=(-\Delta)^{-1}\operatorname{div}(\rho\dot u)+(2\mu+\lambda)\operatorname{div} u.
\end{align}
Multiplying \eqref{eqp3} by $P$ and integrating by parts, we deduce from H{\"o}lder's inequality that
\begin{align}\label{xz2}
\|P\|_{L^{2}}^{2}&\leq C\left\|(-\Delta)^{-1} \operatorname{div}(\rho \dot{u})\right\|_{L^{4 \gamma}}\|P\|_{L^{4 \gamma/(4 \gamma-1)} }+C\|\nabla u\|_{L^{2}}\|P\|_{L^{2}}.
\end{align}
Owing to $1<\frac{4\gamma}{2\gamma+1}<2$, one obtains from Lemma \ref{hls} and Remark \ref{br}  that
\begin{align}\label{eqp4}
\bigl\|(-\Delta)^{-1}\operatorname{div}(\rho\dot u)\bigr\|_{L^{4\gamma}}
\le C\|\rho\dot u\|_{L^{\frac{4\gamma}{2\gamma+1}}}
\le C\|\rho\|_{L^{2\gamma}}^{1/2}\|\sqrt\rho \dot u\|_{L^{2}},
\end{align}
while the interpolation inequality and \eqref{eqma} imply
\begin{align}\label{eqp5}
\|P\|_{L^{\frac{4\gamma}{4\gamma-1}}}=A\|\rho\|_{L^{\frac{4\gamma^2}{4\gamma-1}}}^\gamma
\le C m_0^{1/2}\|\rho\|_{L^{2\gamma}}^{\gamma-\frac12}.
\end{align}
Hence, substituting \eqref{eqp4} and \eqref{eqp5} into \eqref{xz2} yields
\begin{align*}
\|P\|_{L^{2}}^{2}&\le C A m_0^{1/2}\|\rho\|_{L^{2\gamma}}^{\gamma}
\|\sqrt\rho \dot u\|_{L^{2}}+C\|\nabla u\|_{L^{2}}\|P\|_{L^{2}}\notag\\
&=C m_0^{1/2}\|P\|_{L^{2}}\|\sqrt\rho \dot u\|_{L^{2}}+C\|\nabla u\|_{L^{2}}\|P\|_{L^{2}}\notag\\
&\leq \frac12\|P\|_{L^{2}}^{2}+ C m_0\|\sqrt\rho \dot u\|_{L^{2}}^{2}+C\|\nabla u\|_{L^{2}}^{2},
\end{align*}
proving \eqref{eqp1}. Integrating \eqref{eqp1} in time and using \eqref{eqhe}, we obtain
\eqref{eqp2}.
\end{proof}

\begin{lemma}
	It holds that
	\begin{align}
		\sup_{0\le t\le T}D(t)
		+\int_0^T\left(\|\sqrt\rho \dot u\|_{L^{2}}^{2}+\|P\|_{L^{3}}^3\right)d t
		\le CD_{0}+C\bar\rho\int_0^TD(t)^2 d t.
		\label{eq40}
	\end{align}
\end{lemma}

\begin{proof}
	It follows from \eqref{eqf1} and \eqref{eq32} that
	\begin{align}
		&\frac{  d}{  dt}
		\left[
		\frac\mu2\|\nabla u\|_{L^{2}}^{2}
		+\frac{\mu+\lambda}{2}\|\operatorname{div} u\|_{L^{2}}^{2}
		-\int P\operatorname{div} u d x
		\right]
		+\|\sqrt\rho \dot u\|_{L^{2}}^{2}
		\notag\\
		&\le C\int P|\nabla u|^2 d x
		+C\|\nabla u\|_{L^{3}}^3
		\notag\\
		&\le C\|P\|_{L^{3}}^3+C\|\nabla u\|_{L^{3}}^3
		\notag\\[2mm]
		&\le C\|P\|_{L^{3}}^3
		+C\bar\rho^{1/2}\|\sqrt\rho \dot u\|_{L^{2}}D(t).
		\label{eq41}
	\end{align}
	
The continuity equation and the pressure law imply
\begin{align}\label{eqpre}
P_t+u\cdot\nabla P+\gamma P\operatorname{div} u=0.
\end{align}
Multiplying \eqref{eqpre} by $2P$ and integrating by parts, we find
\begin{align*}
\frac{d}{dt}\|P\|_{L^{2}}^{2}+\frac{2\gamma-1}{2\mu+\lambda}\|P\|_{L^{3}}^3
=-\frac{2\gamma-1}{2\mu+\lambda}\int P^2F d x,
\end{align*}
which along with Young's inequality and \eqref{eq31} yields
\begin{align}\label{eq42}
\frac{d}{dt}\|P\|_{L^{2}}^{2}+\frac{2\gamma-1}{2(2\mu+\lambda)}\|P\|_{L^{3}}^3
\le C\|F\|_{L^{3}}^3\le C\bar\rho^{1/2}\|\sqrt\rho \dot u\|_{L^{2}}D(t).
\end{align}

Define
\begin{align}\label{eq43}
\mathcal E_1(t)\triangleq\frac\mu2\|\nabla u\|_{L^{2}}^{2}
+\frac{\mu+\lambda}{2}\|\operatorname{div} u\|_{L^{2}}^{2}-\int P\operatorname{div} u d x
+K_P\|P\|_{L^{2}}^{2},
\end{align}
where $K_P>0$ is a constant depending only on $\mu$, $\lambda$, and $\gamma$. Choose $K_P$ sufficiently large so that \eqref{eq42} multiplied by $K_P$ can absorb the pressure term on the right-hand side of \eqref{eq41}, and hence
\begin{align}\label{eq44}
\frac{d}{dt}\mathcal E_1(t)+\|\sqrt\rho \dot u\|_{L^{2}}^{2}+\|P\|_{L^{3}}^3
\le C\bar\rho^{1/2}\|\sqrt\rho \dot u\|_{L^{2}}D(t)\le\frac12\|\sqrt\rho \dot u\|_{L^{2}}^{2}
+C\bar\rho D(t)^2.
\end{align}	
Since
\begin{align*}
\left|\int P\operatorname{div} u dx\right|\le\frac\mu4\|\nabla u\|_{L^{2}}^{2}+C\|P\|_{L^{2}}^{2},
\end{align*}
increasing $K_P$, if necessary, gives
\begin{align}\label{eq45}
\mathcal E_1(t)\sim  D(t).
\end{align}
Integrating \eqref{eq44} and using \eqref{eq45}, we obtain \eqref{eq40}.
\end{proof}

\subsection{Scaling-invariant bootstrap estimates}
The lower-order control functional is
\begin{align*}
\begin{split}
\mathcal S(T)\triangleq \bar\rho E_0+\bar\rho m_0\left\{\sup_{0\le t\le T}D(t)
+\int_0^T\left(\|\sqrt\rho \dot u\|_{L^{2}}^{2}+\|P\|_{L^{3}}^3\right) d t\right\}.
\end{split}
\end{align*}

\begin{lemma}[Bootstrap estimate]\label{bootstrap}
There exist constants $K\ge2$ and $\varepsilon_*>0$ depending only on $\mu, \lambda$, and $\gamma$ such that, if
\begin{align}\label{cpp1}
0\le\rho\le2\bar\rho,\quad \mathcal S(T)\le2K \mathcal S_0,
\end{align}
then
\begin{gather}\label{cpp2}
\mathcal S(T)\le K \mathcal S_0,
\\ \label{cp50}
\int_0^TD(t)d t\le C\bigl(E_0+m_0D_{0}\bigr)
\end{gather}
provided $S_0\le\varepsilon_*$.
\end{lemma}
\begin{proof}\label{cp51}
It follows from \eqref{eqp2} and \eqref{eq40} that
\begin{align}
\sup_{0\le t\le T}D(t)+\int_0^T\left(\|\sqrt\rho \dot u\|_{L^{2}}^{2}+\|P\|_{L^{3}}^3\right) dt
\le CD_{0}+C\bar\rho\sup_{0\le t\le T}D(t)\left[E_0+m_0\int_0^T\|\sqrt\rho \dot u\|_{L^{2}}^{2} d t\right].
\end{align}
Multiplying \eqref{cp51} by $\bar\rho m_0$ and adding $\bar\rho E_0$ gives
\begin{align}\label{cp52}
\mathcal S(T) \le C_1 \mathcal S_0+C_2\bar\rho m_0\sup_{0\le t\le T}D(t)\left(\bar\rho E_0+
\bar\rho m_0\int_0^T\|\sqrt\rho \dot u\|_{L^{2}}^{2}d t\right) \le C_1 \mathcal S_0+C_2 \mathcal S(T)^2,
\end{align}
where $C_1\ge1$ and $C_2>0$ depend only on $\mu, \lambda$, and $\gamma$. Set $K\triangleq2C_1$ and $\varepsilon_*\triangleq\frac1{8C_2K}$, then one deduces from \eqref{cpp1} that
\begin{align*}
C_2 \mathcal S(T)^2\le2C_2K\mathcal S_0 \mathcal S(T)\le\frac14\mathcal S(T)
\end{align*}
provided $S_0\le\varepsilon_*$, which combined with \eqref{cp52} implies
\begin{align*}
\mathcal S(T)\le\frac43C_1 \mathcal S_0\le K \mathcal S_0,
\end{align*}
proving \eqref{cpp2}.
	
Moreover, one gets from \eqref{eqp2} and \eqref{cpp2} that
\begin{align*}
\int_0^TD(t) d t\le CE_0+Cm_0\int_0^T\|\sqrt\rho \dot u\|_{L^{2}}^{2} d t
\le CE_0+C\bar\rho^{-1}\mathcal S(T) \le C\bigl(E_0+m_0D_{0}\bigr),
\end{align*}
as asserted.
\end{proof}

\begin{remark}[Scaling and the mass factor]\label{rem}
The factor $m_0$ in $\mathcal S_0$ is determined by the pressure estimate \eqref{eqp1}. After time integration, it multiplies the material-derivative dissipation in \eqref{eqp2}; this is precisely the factor needed in \eqref{cp51}.
\end{remark}

We continue the lower-order estimates by deriving a uniform upper bound for the density. The argument combines time-weighted energy estimates, mass localization through a first spatial moment, and control of the effective viscous flux along particle trajectories.

Since
\begin{align*}
|x|\leq \bar{x}\leq \bar{x}^a,\quad x\in\mathbb R^2,
\end{align*}
the condition $\bar{x}^a\rho_0\in L^1$ implies
\begin{align}\label{kk}
\int\rho_0(x)|x|dx\leq\int\rho_0(x)\bar{x}^adx<\infty.
\end{align}

We recall the associated initial quantities
\begin{align*}
L_0 =\left(\frac{m_0}{\bar\rho}\right)^{1/2},
\quad\mathcal M_0 =\frac1{m_0L_0}\int\rho_0(x)\sqrt{L_0^2+|x|^2} d x,
\quad \Pi_0= A m_0\bar\rho^\gamma.
\end{align*}
Thus $1\le\mathcal M_0<\infty$. The smallness quantity used in this section is
\begin{align*}
\mathcal S_0^{\mathrm{vac}}=\mathcal M_0^{1/3}\bigl(1+\Pi_0^{-2/3}\bigr)\mathcal S_0.
\end{align*}
In particular, $\mathcal S_0\le\mathcal S_0^{\mathrm{vac}}$. Under the exact scaling \eqref{2sca} (see \ref{dd}),
\begin{align*}
\mathcal M_0^\ell=\mathcal M_0,
\quad \Pi_0^\ell=\Pi_0,
\quad \left( \mathcal S_0^{\mathrm{vac}}\right)^\ell=\mathcal S_0^{\mathrm{vac}}.
\end{align*}

\begin{proposition}[Bootstrap estimate]\label{prop}
Let $K\ge2$ be the constant in Lemma \ref{bootstrap}. There exists $0<\varepsilon_0\in\min\{1, \varepsilon_*\}$ depending only on $\mu, \lambda$, and $\gamma$ such that, if the bootstrap assumption \eqref{cpp1} holds and
\begin{align}\label{bao}
\mathcal S_0^{\mathrm{vac}}\le\varepsilon_0,
\end{align}
then
\begin{align*}
\sup_{0\le t\le T}\|\rho(t)\|_\infty\le\frac74\bar\rho.
\end{align*}
\end{proposition}
\begin{proof}
The proof is given in subsection \ref{ppp}.
\end{proof}

Since $\mathcal S_0\le\mathcal S_0^{\mathrm{vac}}$, one derives from Lemma \ref{bootstrap} that
\begin{align}\label{de}
\mathcal S(T)\le 2K\mathcal S_0,\quad \int_0^T D(t) dt
\le C\bigl(E_0+m_0D_{0}\bigr),\quad \bar\rho\int_0^TD(t) dt\le C\mathcal S_0.
\end{align}

\begin{lemma}[Two time-weighted estimates]
Under \eqref{bao}, it holds that
\begin{gather}\label{60}
\sup_{0\le t\le T}\big[(m_0+t)D(t)\big]+\int_0^T(m_0+t)
\big(\|\sqrt\rho \dot u\|_{L^{2}}^{2}+\|P\|_{L^{3}}^3\big) d t\le C\bigl(E_0+m_0D_{0}\bigr),
\\
\sup_{0<t\le T}\big[t(m_0+t)\big(\|\sqrt\rho \dot u\|_{L^{2}}^{2}+\|P\|_{L^{3}}^3\big)\big]
+\int_0^Tt(m_0+t)\big(\|\nabla\dot u\|_{L^{2}}^{2}+\|P\|_{L^{4}}^4\big) d t\le C\bigl(E_0+m_0D_{0}\bigr).
\label{61}
\\
D(t)\le\frac{C\mathcal S_0}{\bar\rho(m_0+t)},\quad
\|\sqrt\rho \dot u(t)\|_{L^{2}}^{2}\le\frac{C\mathcal S_0}{\bar\rho\,t(m_0+t)},\quad \text{for}\quad 0<t\le T.
\label{62}
\end{gather}
\end{lemma}
\begin{proof}
Multiplying \eqref{eq44} by $(m_0+t)$ and using \eqref{eq45}, one has that
\begin{align}\label{63}
\frac{  d}{  dt}\bigl[(m_0+t)\mathcal E_1(t)\bigr]+c(m_0+t)\left(\|\sqrt\rho \dot u\|_{L^{2}}^{2}+\|P\|_{L^{3}}^3\right)\le\mathcal E_1(t)+C\bar\rho D(t)(m_0+t)\mathcal E_1(t).
\end{align}
By \eqref{de},
\begin{align*}
m_0\mathcal E_1(0)+\int_0^T\mathcal E_1(t) d t
\le C\bigl(E_0+m_0D_{0}\bigr),\quad \bar\rho\int_0^TD(t) d t\le C\mathcal S_0\leq C,
\end{align*}
which combined with \eqref{63} and Gronwall's inequality yields \eqref{60}.
	
It follows from \eqref{c19}, \eqref{eq34}, and Young's inequality that
\begin{align}\label{64}
&\frac12\frac{d}{dt}\|\sqrt\rho \dot u\|_{L^{2}}^{2}+\mu\|\nabla\dot u\|_{L^{2}}^{2}
+(\mu+\lambda)\|\operatorname{div}\dot u\|_{L^{2}}^{2}\\
&\le C\int\left(|\nabla u|^2+P|\nabla u|\right)|\nabla\dot u| dx\\
&\le \frac{\mu}{2}\|\nabla\dot u\|_{L^{2}}^{2}+C\|\nabla u\|_{L^{4}}^4+C\|P\|_{L^{4}}^4\\
&\le \frac{\mu}{2}\|\nabla\dot u\|_{L^{2}}^{2}+C\bar\rho D(t)\|\sqrt\rho \dot u\|_{L^{2}}^{2}+C\|P\|_{L^{4}}^4.
\end{align}
Multiplying \eqref{eqpre} by $3P^2$ and integration by parts, we obtain from Young's inequality and \eqref{eqe33} that
\begin{align}\label{65}
\frac{d}{dt}\|P\|_{L^{3}}^3+\|P\|_{L^{4}}^4\le C\|F\|_{L^{4}}^4
\le C\bar\rho D(t)\|\sqrt\rho \dot u\|_{L^{2}}^{2}.
\end{align}
Choose a fixed sufficiently large $K_3>0$ and set
\begin{align*}
\mathcal E_2(t)\triangleq\|\sqrt\rho \dot u\|_{L^{2}}^{2}+K_3\|P\|_{L^{3}}^3.
\end{align*}
Adding \eqref{65} multiplied by $K_3$ to \eqref{64}, one sees that
\begin{align*}
\frac {d}{dt}\mathcal E_2 (t)+\|\nabla\dot u\|_{L^{2}}^{2}+\|P\|_{L^{4}}^4
\le C\bar\rho D(t)\mathcal E_2(t).
\end{align*}
Multiplication by $t(m_0+t)$ yields
\begin{align*}
\frac{d}{dt}\bigl[t(m_0+t)\mathcal E_2(t)\bigr]+ct(m_0+t)\left(\|\nabla\dot u\|_{L^{2}}^{2}+\|P\|_{L^{4}}^4\right)
\le C\bar\rho D(t)t(m_0+t)\mathcal E_2(t)+2(m_0+t)\mathcal E_2(t),
\end{align*}
which combined with \eqref{de}, \eqref{60}, and Gronwall's inequality implies \eqref{61}. The desired \eqref{62} follows from \eqref{60} and \eqref{61}.
\end{proof}

\begin{lemma}[Mass localization]
It holds that, for $0\le t\le T$,
\begin{gather}\label{70}
\mathcal M(t)\le C\mathcal M_0\left(1+\frac{t}{m_0}\right),\\
\int_{B_{4L_0\mathcal M(t)}}\rho(x,t) dx\ge\frac34m_0.\label{71}
\end{gather}
\end{lemma}
\begin{proof}
Since
\begin{align*}
\nabla_x\sqrt{L_0^2+|x|^2}=\frac{x}{\sqrt{L_0^2+|x|^2}},
\end{align*}
we have
\begin{align*}
\left|\nabla_x\sqrt{L_0^2+|x|^2}\right|=\frac{|x|}{\sqrt{L_0^2+|x|^2}}\leq 1,
\end{align*}
which together with the continuity equation and \eqref{eqhe} gives
\begin{align*}
\frac{d}{dt}\int\rho\sqrt{L_0^2+|x|^2} d x\le\int\rho|u| d x
\le\sqrt{m_0}\,\|\sqrt\rho u\|_{L^{2}}
\le\sqrt{2m_0E_0}.
\end{align*}
Dividing the above estimate by $m_0L_0$ and integrating with respect to time, we obtain
\begin{align*}
\mathcal M(t)\leq \mathcal M_0+\sqrt{2\bar\rho E_0}\,\frac{t}{m_0},
\end{align*}
proving \eqref{70} due to $\bar\rho E_0\leq \mathcal S_0\leq1$ and $\mathcal M_0\geq1$.

Moreover, one sees that
\begin{align*}
\int_{|x|>4L_0\mathcal M(t)}\rho d x\le\frac{1}{4L_0\mathcal M(t)}\int\rho|x| d x\le\frac14m_0,
\end{align*}
which implies \eqref{71}.
\end{proof}

\begin{lemma}[Weighted integrability with explicit parameters]\label{zuilem}
Let $2\le r<\infty$. If $\nabla v\in L^2$ and $\sqrt\rho v\in L^2$, then
\begin{align*}
\|\rho v\|_{L^r}\le C_r\mathcal M_0^4\left(1+\frac{t}{m_0}\right)^4\bar\rho^{1-1/r}m_0^{1/r}
\big(m_0^{-1/2}\|\sqrt\rho v\|_{L^{2}}+\|\nabla v\|_{L^{2}}\big).
\end{align*}
\end{lemma}
\begin{proof}
\emph{Step 1: Control of the localization of mass.}
Introduce
\begin{align*}
\bar x_{L_0}\triangleq\left(e+\frac{|x|^2}{L_0^2}\right)^{1/2}\log^2\left(e+\frac{|x|^2}{L_0^2}\right).
\end{align*}
Since $\log z\le Cz^{1/4}$ for $z\ge e$, we have
\begin{align*}
(e+r^2)^{1/4}\log(e+r^2)\le C(e+r^2)^{1/2}\le C\sqrt{1+r^2},\quad r\ge0.
\end{align*}
Consequently,
\begin{align}\label{80}
\frac1{m_0}\int\rho\bar x_{L_0}^{1/2} d x\le\frac{C}{m_0}\int\rho\sqrt{1+\frac{|x|^2}{L_0^2}} d x
=\frac{C}{m_0L_0}\int\rho\sqrt{L_0^2+|x|^2} d x=C\mathcal M(t).
\end{align}
Also, it follows from $\int  \rho(x, t)|x| \mathrm{d} x \leq m_{0} L_{0} \mathcal{M}(t)$ that
\begin{align*}
\int_{|x|>4L_0\mathcal M(t)}\rho d x\le\frac1{4L_0\mathcal M(t)}\int\rho|x| d x
\le\frac{m_0L_0\mathcal M(t)}{4L_0\mathcal M(t)}=\frac14m_0,
\end{align*}
which leads to
\begin{align}\label{81}
\int_{B_{4L_0\mathcal M(t)}}\rho d x\ge\frac34m_0.
\end{align}
	
\emph{Step 2: A local $L^2$ estimate from a ball by mass.}
Let $R$ be any radius satisfying
\begin{align}\label{82}
R\ge L_0,\quad \int_{B_R}\rho d x\ge\frac34m_0.
\end{align}
Denote the spatial average of $v$ over $B_R$ by $v_{B_R}\triangleq\frac1{|B_R|}\int_{B_R}v d x$, then one deduces from
\begin{align*}
\left(\int_{B_R}\rho d x\right)v_{B_R}=\int_{B_R}\rho v d x+\int_{B_R}\rho(v_{B_R}-v) d x
\end{align*}
and \eqref{82} that
\begin{align*}
\frac34m_0|v_{B_R}|&\le\left(\int_{B_R}\rho d x\right)^{1/2}
\|\sqrt\rho v\|_{L^2(B_R)}+\|\rho\|_{L^2(B_R)}\|v-v_{B_R}\|_{L^2(B_R)}
\\&\le C\bar\rho^{1/2}R\|\sqrt\rho v\|_{L^{2}}+C\bar\rho R^2\|\nabla v\|_{L^{2}},
\end{align*}
due to Cauchy--Schwarz inequality, $\rho\le2\bar\rho$, $|B_R|=\pi R^2$, and Poincar\'e's inequality
\begin{align*}
\|v-v_{B_R}\|_{L^2(B_R)}\le CR\|\nabla v\|_{L^2(B_R)}.
\end{align*}
Since $\bar\rho=m_0/L_0^2$, division by $m_0$ gives
\begin{align*}
|v_{B_R}|\le C\frac{R}{L_0}m_0^{-1/2}\|\sqrt\rho v\|_{L^{2}}+C\left(\frac{R}{L_0}\right)^2\|\nabla v\|_{L^{2}}.
\end{align*}
Application of Poincar\'e's inequality and $R/L_0\ge1$ implies
\begin{align}\label{83}
L_0^{-1}\|v\|_{L^2(B_{L_0})}\le C\frac{R}{L_0}\|\nabla v\|_{L^{2}}+C\frac{R}{L_0}|v_{B_R}|
\le C\left(\frac{R}{L_0}\right)^3\big(m_0^{-1/2}\|\sqrt\rho v\|_{L^{2}}+\|\nabla v\|_{L^{2}}\big).
\end{align}
	
\emph{Step 3: The weighted Sobolev estimate.}
Fix $2\le q<\infty$ and choose $\theta=q+2$. For $z\ge e$,
\begin{align}\label{xz4}
z^{-2}(\log z)^{-8}\le c_q z^{-1}(\log z)^{-(q+2)}.
\end{align}
Taking $z=e+|y|^2$, multiplying \eqref{xz4} by $|f(y)|^q$, and integrating, we obtain from Lemma \ref{lious} that
\begin{align}\label{84}
&\left(\int_{\R^2}|f(y)|^q(e+|y|^2)^{-2}\log^{-8}(e+|y|^2)\, d y\right)^{1/q}
\notag\\& \le c_q^{1/q}\left(\int_{\R^2}\frac{|f(y)|^q}{e+|y|^2}\log^{-(q+2)}(e+|y|^2)\, d y\right)^{1/q}
\notag\\& \le C_q\left(\|f\|_{L^2(B_1)}+\|\nabla_y f\|_{L^2(\R^2)}\right).
\end{align}
Set $x=L_0y$ and $f(y)=v(L_0y)$, then one has
\begin{align}\label{85}
\|f\|_{L^2(B_1)}^2=\int_{|y|<1}|v(L_0y)|^2\, d y=L_0^{-2}\int_{|x|<L_0}|v(x)|^2\, d x.
\end{align}
The chain rule gives
\begin{align*}
 \partial_{y_j}f(y)=L_0(\partial_{x_j}v)(L_0y),\quad j=1,2.
\end{align*}
Consequently,
\begin{align}
\|\nabla_y f\|_{L^2(\R^2)}^2=L_0^2\int_{\R^2}|\nabla v(L_0y)|^2\,d y
=L_0^2L_0^{-2}\int_{\R^2}|\nabla v(x)|^2\,d x=\|\nabla v\|_{L^{2}}^{2}.
\label{86}
\end{align}
Noting that $\bar x_{L_0}(L_0y)^{-4}=(e+|y|^2)^{-2}\log^{-8}(e+|y|^2)$, and hence
\begin{align*}
\frac1{L_0^2}\int_{\R^2}|v(x)|^q\bar x_{L_0}(x)^{-4}dx
=\int_{\R^2}|v(L_0y)|^q\bar x_{L_0}(L_0y)^{-4}dy
=\int_{\R^2}|f(y)|^q(e+|y|^2)^{-2}\log^{-8}(e+|y|^2)dy,
\end{align*}
which together with \eqref{84}--\eqref{86} and \eqref{83} shows
\begin{align}\label{88}
\left(\frac1{L_0^2}\int|v|^q\bar x_{L_0}^{-4}\,d x\right)^{1/q}
&\le C_q\left(L_0^{-1}\|v\|_{L^2(B_{L_0})}+\|\nabla v\|_{L^{2}}\right)
\notag\\&\le C_q\left[\left(\frac{R}{L_0}\right)^3\left(
m_0^{-1/2}\|\sqrt\rho v\|_{L^{2}}+\|\nabla v\|_{L^{2}}\right)
+\|\nabla v\|_{L^{2}}\right]
\notag\\&\le C_q\left(\frac{R}{L_0}\right)^3
\big(m_0^{-1/2}\|\sqrt\rho v\|_{L^{2}}+\|\nabla v\|_{L^{2}}\big).
\end{align}

\emph{Step 4: Recovering the density-weighted $L^r$ norm.}
We obtain from H\"older's inequality that	
\begin{align*}
\int\rho|v|^r d x=\int\bigl(\rho\bar x_{L_0}^{1/2}\bigr)^{8/9}\rho^{1/9}|v|^r\bar x_{L_0}^{-4/9} d x
\le C\bar\rho^{1/9}\left(\int\rho\bar x_{L_0}^{1/2} d x\right)^{8/9}\left(\int|v|^{9r}\bar x_{L_0}^{-4} d x\right)^{1/9}.
\end{align*}
After dividing by $m_0$ and taking the $r$-th root, one infers from \eqref{88} with $q=9r$ and Young's inequality that
\begin{align}\label{89}
\left(\frac1{m_0}\int\rho|v|^r d x\right)^{1/r}
& \le C\left(\frac1{m_0}\int\rho\bar x_{L_0}^{1/2} d x\right)^{8/(9r)}
\left(\frac1{L_0^2}\int|v|^{9r}\bar x_{L_0}^{-4} d x\right)^{1/(9r)} \notag \\
&\le C_r\left(\frac{R}{L_0}\right)^3
\left(1+\frac1{m_0}\int\rho\bar x_{L_0}^{1/2} d x\right)\big(m_0^{-1/2}\|\sqrt\rho v\|_{L^{2}}+\|\nabla v\|_{L^{2}}\big).
\end{align}

\emph{Step 5: Choice of the radius and the final estimate.}
Choose $R=4L_0\mathcal M(t)$, which satisfies \eqref{82} by \eqref{81}. By \eqref{80}
and $\mathcal M(t)\ge1$,
\begin{align*}
\left(\frac{R}{L_0}\right)^3\left(1+\frac1{m_0}\int\rho\bar x_{L_0}^{1/2}\, d x\right)
\le C\mathcal M(t)^3\bigl(1+\mathcal M(t)\bigr)\le C\mathcal M(t)^4.
\end{align*}
Multiplying \eqref{89} by $m_0^{1/r}$, one has that
\begin{align*}
\left(\int\rho|v|^r d x\right)^{1/r}\le C_r m_0^{1/r}\mathcal M(t)^4
\big(m_0^{-1/2}\|\sqrt\rho v\|_{L^{2}}+\|\nabla v\|_{L^{2}}\big),
\end{align*}
which combined with \eqref{70} implies
\begin{align*}
\|\rho v\|_{L^r}&\le(2\bar\rho)^{1-1/r}\left( \int\rho|v|^r d x\right)^{1/r}\\
&\leq C_r (\bar\rho)^{1-1/r} m_0^{1/r}\mathcal M(t)^4
\big(m_0^{-1/2}\|\sqrt\rho v\|_{L^{2}}+\|\nabla v\|_{L^{2}}\big)\\
&\le C_r\mathcal M_0^4\bigg(1+\frac{t}{m_0}\bigg)^4
\bar\rho^{1-1/r}m_0^{1/r}
\big(m_0^{-1/2}\|\sqrt\rho v\|_{L^{2}}+\|\nabla v\|_{L^{2}}\big). \tag*{\qedhere}
\end{align*}
\end{proof}

\subsection{Short-time and large-time effective viscous flux estimates}

\begin{lemma}[Short-time and large-time flux bounds]\label{flux}
Under \eqref{bao}, it holds that
\begin{gather}
\int_0^{\min\{m_0,T\}}\|F(t)\|_\infty d t\le C\mathcal M_0^{1/6}\mathcal S_0^{1/2},\label{90} \\
m_0^{1/2}\int_{m_0}^T\|F(t)\|_\infty^{3/2} d t\le C\mathcal M_0^{1/4}\mathcal S_0^{3/4},\quad \text{for}\quad T>m_0.\label{91}
\end{gather}
\end{lemma}
\begin{proof}
\emph{Step 1: An $L^{48}$ estimate for $F$.}
It follows from \eqref{L2}, \eqref{t2}, and \eqref{62} that
\begin{gather}
\|F(t)\|_{L^{2}}\le\frac{C\mathcal S_0^{1/2}}{\bar\rho^{1/2}(m_0+t)^{1/2}}
=\frac{C\mathcal S_0^{1/2}L_0}{m_0}\left(1+\frac{t}{m_0}\right)^{-1/2},
\label{92}\\
\|\nabla F(t)\|_{L^{2}}\le C\bar\rho^{1/2}
\frac{\mathcal S_0^{1/2}}{\bar\rho^{1/2}t^{1/2}(m_0+t)^{1/2}}=\frac{C\mathcal S_0^{1/2}}{m_0}
\left(\frac{t}{m_0}\right)^{-1/2}\left(1+\frac{t}{m_0}\right)^{-1/2},
\label{93}
\end{gather}
which combined with Gagliardo--Nirenberg inequality indicates
\begin{align}
\|F(t)\|_{L^{48}}&\le C\|F(t)\|_{L^{2}}^{1/24}\|\nabla F(t)\|_{L^{2}}^{23/24}\notag\\
&\le C\bigg[\frac{\mathcal S_0^{1/2}L_0}{m_0}
\left(1+\frac{t}{m_0}\right)^{-1/2}\bigg]^{1/24}
\bigg[\frac{\mathcal S_0^{1/2}}{m_0}
		\left(\frac{t}{m_0}\right)^{-1/2}
		\left(1+\frac{t}{m_0}\right)^{-1/2}\bigg]^{23/24}
		\notag\\
		&=\frac{C\mathcal S_0^{1/2}L_0^{1/24}}{m_0}
		\left(\frac{t}{m_0}\right)^{-23/48}
		\left(1+\frac{t}{m_0}\right)^{-1/2}.
		\label{94}
	\end{align}
	
\emph{Step 2: An $L^{48}$ estimate for $\nabla F$.}
Applying Lemma \ref{zuilem} and \eqref{eqer} with $r=48$ and $v=\dot u$ to get
	\begin{align}
		\|\nabla F(t)\|_{48}
		\le C\|\rho\dot u(t)\|_{48}
		&\le C\mathcal M_0^4\left(1+\frac{t}{m_0}\right)^4
		\bar\rho^{47/48}m_0^{1/48}
		\big(m_0^{-1/2}\|\sqrt\rho \dot u(t)\|_{L^{2}}
		+\|\nabla\dot u(t)\|_{L^{2}}\big)\notag\\
		&\le\frac{C\mathcal M_0^4}{m_0L_0^{23/24}}
		\left(1+\frac{t}{m_0}\right)^4
		\left[\sqrt{\bar\rho}\,m_0
		\left(\|\sqrt\rho \dot u(t)\|_{L^{2}}
		+\sqrt{m_0}\|\nabla\dot u(t)\|_{L^{2}}\right)\right]
		\label{95},
	\end{align}
due to
\begin{align*}
\bar\rho^{47/48}m_0^{1/48}\big(m_0^{-1/2}\|\sqrt\rho \dot u\|_{L^{2}}+\|\nabla\dot u\|_{L^{2}}\big)
&=\bar\rho^{47/48}m_0^{-23/48}\big(\|\sqrt\rho \dot u\|_{L^{2}}+\sqrt{m_0}\|\nabla\dot u\|_{L^{2}}\big)
\notag\\&=\frac1{m_0L_0^{23/24}}\left[\sqrt{\bar\rho}\,m_0
\left(\|\sqrt\rho \dot u\|_{L^{2}}+\sqrt{m_0}\|\nabla\dot u\|_{L^{2}}\right)\right].
\end{align*}
	
\emph{Step 3: Interpolation into $L^\infty$.} It follows from \eqref{g2} that
\begin{align*}
\|F\|_\infty\le C\|F\|_{48}^{23/24}\|\nabla F\|_{48}^{1/24},
\end{align*}
which together with \eqref{94} and \eqref{95} leads to
\begin{align}\label{wuq}
\|F(t)\|_{L^\infty}\le\frac{C\mathcal M_0^{1/6}\mathcal S_0^{23/48}}{m_0}
\left(\frac{t}{m_0}\right)^{-\frac{529}{1152}}\left(1+\frac{t}{m_0}\right)^{-\frac{5}{16}}
\left[\sqrt{\bar\rho}\,m_0\left(\|\sqrt\rho \dot u(t)\|_{L^{2}}+\sqrt{m_0}\|\nabla\dot u(t)\|_{L^{2}}\right)\right]^{\frac{1}{24}}.
\end{align}
	
\emph{Step 4: The weighted square integral.} Define $\sigma_{m_0}(t)\triangleq\min\big\{1,\frac{t}{m_0}\big\}$,
then one has that
\begin{align}\label{96}
&\int_0^T\sigma_{m_0}(t)\left[\sqrt{\bar\rho}\,m_0
\left(\|\sqrt\rho \dot u\|_{L^{2}}+\sqrt{m_0}\|\nabla\dot u\|_{L^{2}}\right)\right]^2
\frac{ d t}{m_0}\notag\\
& \le2\bar\rho m_0\int_0^T\sigma_{m_0}(t)\|\sqrt\rho \dot u\|_{L^{2}}^{2} d t
+2\bar\rho m_0^2\int_0^T\sigma_{m_0}(t)\|\nabla\dot u\|_{L^{2}}^{2} d t\notag\\
& \le2\bar\rho m_0\int_0^T\|\sqrt\rho \dot u\|_{L^{2}}^{2}\, d t
+2\bar\rho\int_0^Tt(m_0+t)\|\nabla\dot u\|_{L^{2}}^{2}\, d t
\notag\\&\le C\mathcal S_0,
\end{align}
where we have used \eqref{de}, \eqref{61}, and
\begin{align*}
m_0^2\sigma_{m_0}(t)=\begin{cases}	
m_0t,&0\le t\le m_0,\\
m_0^2,&t\ge m_0,
\end{cases}
\quad
m_0^2\sigma_{m_0}(t)\le t(m_0+t).
\end{align*}
Changing the integration variable to $s=t/m_0$ in \eqref{96} gives
\begin{align}\label{97}
&\int_0^{T/m_0}\min\{1,s\}\left[\sqrt{\bar\rho}\,m_0
\left(\|(\sqrt\rho \dot u)(m_0s)\|_{L^{2}}+\sqrt{m_0}\|\nabla\dot u(m_0s)\|_{L^{2}}\right)\right]^2\, d s
\le C\mathcal S_0.
\end{align}
	
\emph{Step 5: The short-time estimate.} Set $s_*\triangleq\min\{1,T/m_0\}$. From
\eqref{wuq}, H\"older's inequality, and \eqref{97}, one sees that
\begin{align*}
&\int_0^{\min\{m_0,T\}}\|F(t)\|_{L^\infty} dt\notag\\
&\le C\mathcal M_0^{\frac16}\mathcal S_0^{\frac{23}{48}}
\int_0^{s_*}s^{-\frac{529}{1152}}(1+s)^{-\frac{5}{16}}
\left[\sqrt{\bar\rho}\,m_0\left(\|(\sqrt\rho \dot u)(m_0s)\|_{L^{2}}
+\sqrt{m_0}\|\nabla\dot u(m_0s)\|_{L^{2}}\right)\right]^{\frac{1}{24}} d s\\
&\le C\mathcal M_0^{\frac{1}{6}}\mathcal S_0^{\frac{23}{48}}
\left(\int_0^{s_*}s^{-553/1128}\, d s\right)^{\frac{47}{48}}
\left\{\int_0^{s_*}s\left[\sqrt{\bar\rho}\,m_0\left(\|(\sqrt\rho \dot u)(m_0s)\|_{L^{2}}
+\sqrt{m_0}\|\nabla\dot u(m_0s)\|_{L^{2}}\right)\right]^2\, d s
\right\}^{\frac{1}{48}}\notag\\
&\le C\mathcal M_0^{1/6}\mathcal S_0^{1/2},
\end{align*}
proving \eqref{90}.
	
\emph{Step 6: The large-time estimate.} Suppose $T>m_0$. For $s=t/m_0\ge1$, it follows from \eqref{wuq} that
\begin{align*}
\|F(m_0s)\|_{L^\infty}^{3/2}\le\frac{C\mathcal M_0^{1/4}\mathcal S_0^{23/32}}{m_0^{3/2}}
s^{-889/768}\left[\sqrt{\bar\rho}\,m_0\left(\|(\sqrt\rho \dot u)(m_0s)\|_{L^{2}}
+\sqrt{m_0}\|\nabla\dot u(m_0s)\|_{L^{2}}\right)\right]^{1/16},
\end{align*}
which along with H\"older's inequality and \eqref{97} yields
\begin{align*}
&m_0^{1/2}\int_{m_0}^T\|F(t)\|_{L^\infty}^{3/2}\, d t\notag\\
&\le C\mathcal M_0^{1/4}\mathcal S_0^{23/32}\int_1^{T/m_0}s^{-889/768}
\left[\sqrt{\bar\rho}\,m_0\left(\|(\sqrt\rho \dot u)(m_0s)\|_{L^{2}}
+\sqrt{m_0}\|\nabla\dot u(m_0s)\|_{L^{2}}\right)\right]^{1/16} d s\\
&\le C\mathcal M_0^{\frac{1}{4}}\mathcal S_0^{\frac{23}{32}}
\bigg(\int_1^{T/m_0}s^{-\frac{889}{744}} d s\bigg)^{\frac{31}{32}}
\bigg\{\int_1^{T/m_0}\left[\sqrt{\bar\rho}\,m_0
\left(\|(\sqrt\rho \dot u)(m_0s)\|_{L^{2}}+\sqrt{m_0}\|\nabla\dot u(m_0s)\|_{L^{2}}\right)\right]^2 d s
\bigg\}^{\frac{1}{32}}\\&\le C\mathcal M_0^{1/4}\mathcal S_0^{3/4},
\end{align*}
as desired \eqref{91}.
\end{proof}

\subsection{Proof of Proposition \ref{prop}}\label{ppp}

Let $X(t)$ be a particle trajectory satisfying $X'(t)=u(X(t),t)$, and set $y(t)=\rho(X(t),t)$. The continuity equation gives
\begin{align}\label{c00}
y'(t)=-\frac{A}{2\mu+\lambda}y(t)^{\gamma+1}-\frac{y(t)}{2\mu+\lambda}F(X(t),t).
\end{align}
Thus, $y^{\prime}(t)=g(y)+b^{\prime}(t)$ holds with
\begin{align*}
g(y)=-\frac{A}{2\mu+\lambda}y^{\gamma+1},
\quad b(t)=-\frac1{2\mu+\lambda}\int_0^t\rho(X(s),s)F(X(s),s) ds.
\end{align*}
	
On the initial interval $[0,\min\{m_0,T\}]$, the first term on the right of \eqref{c00} is nonpositive. Using the assumed bound $\rho\le2\bar\rho$ and \eqref{90}, we obtain
\begin{align*}
\sup_{0\le t\le\min\{m_0,T\}}\|\rho(t)\|_\infty\le\bar\rho
+\frac{2\bar\rho}{2\mu+\lambda}\int_0^{\min\{m_0,T\}}\|F(t)\|_\infty dt
\le\bar\rho+C\bar\rho\mathcal M_0^{1/6}S_0^{1/2}.
\end{align*}
Since
\begin{align*}
\mathcal M_0^{1/6}\mathcal S_0^{1/2}=\big(\mathcal M_0^{1/3}\mathcal S_0\big)^{1/2}
\le\left(\mathcal S_0^{\mathrm{vac}}\right)^{1/2},
\end{align*}
choosing $\varepsilon_0$ in \eqref{bao} sufficiently small gives
\begin{align}\label{cc1}
\sup_{0\le t\le\min\{m_0,T\}}\|\rho(t)\|_{L^{\infty}}\le\frac32\bar\rho.
\end{align}
This already proves the density assertion when $T\le m_0$.
	
Suppose now that $T>m_0$. For any $z\ge0$, Young's inequality $z \leq \delta+C \delta^{-1 / 2} z^{3 / 2}$ gives
\begin{align*}
\frac{2\bar\rho}{2\mu+\lambda}z\le\frac{A\bar\rho^{\gamma+1}}{2\mu+\lambda}+C\bar\rho(A\bar\rho^\gamma)^{-1/2}z^{3/2}.
\end{align*}
Consequently, for $m_0\le t_1<t_2\le T$,
\begin{align*}
|b(t_2)-b(t_1)|&\le\frac{2\bar\rho}{2\mu+\lambda}\int_{t_1}^{t_2}\|F(t)\|_\infty\, d t\notag\\
&\le\frac{A\bar\rho^{\gamma+1}}{2\mu+\lambda}(t_2-t_1)+C\bar\rho(A\bar\rho^\gamma)^{-1/2}
\int_{m_0}^T\|F(t)\|_\infty^{3/2} d t	\notag\\ &\le N_1(t_2-t_1)+N_0,
\end{align*}
where \eqref{91} allows the choices
\begin{align*}
N_1\triangleq\frac{A\bar\rho^{\gamma+1}}{2\mu+\lambda},
\quad N_0\triangleq C\bar\rho\mathcal M_0^{1/4}\Pi_0^{-1/2}\mathcal S_0^{3/4}.
\end{align*}
Note that $g(y)\le-N_1$ for $y\ge\bar\rho$. Then, Zlotnik's inequality and \eqref{cc1} give
\begin{align*}
\sup_{m_0\le t\le T}\|\rho(t)\|_\infty\le\frac32\bar\rho
+C\bar\rho\big(\mathcal M_0^{1/3}\Pi_0^{-2/3}\mathcal S_0\big)^{3/4}
\le\frac32\bar\rho+C\bar\rho \left( \mathcal S_0^{\mathrm{vac}}\right) ^{3/4}
\le\frac74\bar\rho,
\end{align*}
where we have used the smallness of $\varepsilon_0$ in \eqref{bao} such that
\begin{align*}
C\left( \mathcal S_0^{\mathrm{vac}}\right)^{3/4}\leq \frac14.
\end{align*}
We conclude that
\begin{align*}
\sup_{0\le t\le T}\|\rho(t)\|_{L^{\infty}}\le\frac74\bar\rho.
\end{align*}
The bound $\mathcal S(T)\le K\mathcal S_0$ is already contained in \eqref{cpp2}. This completes the proof of Proposition \ref{prop}.\hfill$\square$

\subsection{Proof of Theorem \ref{th2}}

Under the condition \eqref{ini}, a unique local strong solution $(\rho,u)$ to the problem \eqref{1}--\eqref{3} with $\rho_\infty=0$ has been obtained in \cite{LiLiang2014}. Let $T_{\max}>0$ denote its maximal existence time. By Lemma \ref{blowup}, if $T_{\max}<\infty$, then
\begin{align}\label{xz3}
\lim_{T\rightarrow T_{\max}}\|\rho\|_{L^\infty(0,T;L^\infty)}=\infty.
\end{align}
Furthermore, Proposition \ref{prop} together with a standard continuity argument shows that the density bound in \eqref{cpp1} holds on $[0,T]$ for every $T<T_{\max}$ provided \eqref{S1} holds.
This contradicts \eqref{xz3}. Therefore, $T_{\max}=\infty$. \hfill$\square$

\section*{Conflict of interest}
The authors declare that they have no conflict of interest.

\section*{Data availability}
No data was used for the research described in the article.

\section*{Appendix A: the derivation of \eqref{c19}.}
For any sufficiently regular function $f$, one sees
\begin{align*}
\bigl(\partial_t+\partial_k(u^k\,\cdot)\bigr)f=f_t+\partial_k(u^kf)
=f_t+u^k\partial_kf+(\operatorname{div} u)f.
\end{align*}
Applying this operator to
\begin{align*}
\rho\dot u^j=\mu\Delta u^j+(\mu+\lambda)\partial_j\operatorname{div} u-\partial_jP,
\end{align*}
we deduce that
\begin{align}\label{c8}
\partial_t(\rho\dot u^j)+\partial_k(u^k\rho\dot u^j)\notag
&=\mu\big[\Delta u_t^j+\partial_k(u^k\Delta u^j)\big] +(\mu+\lambda)\big[
\partial_j\operatorname{div} u_t+\partial_k(u^k\partial_j\operatorname{div} u)\big]\notag\\&\quad -\partial_jP_t-\partial_k(u^k\partial_jP).
\end{align}
We calculate each terms on the right-hand side of the above identity.

For the Laplacian term, noting that $u_t^j=\dot u^j-u^k\partial_k u^j$, and hence
\begin{align}\label{c10}
\Delta u_t^j=\Delta\dot u^j-\Delta(u^k\partial_k u^j) =\Delta\dot u^j
-\partial_i(\partial_i u^k\partial_k u^j)-\partial_i(u^k\partial_k\partial_i u^j).
\end{align}
While the product rule indicates
\begin{align*}
\partial_k(u^k\Delta u^j)&=\partial_k\big[
\partial_i(u^k\partial_i u^j)-\partial_i u^k\partial_i u^j\big]\notag\\
&=\partial_i\big[\partial_k(u^k\partial_i u^j)\big]
-\partial_k(\partial_i u^k\partial_i u^j)\notag\\&=\partial_i\big[
(\operatorname{div} u)\partial_i u^j+u^k\partial_k\partial_i u^j\big]
-\partial_k(\partial_i u^k\partial_i u^j),
\end{align*}
which together with \eqref{c10} leads to
\begin{align}\label{c12}
\Delta u_t^j+\partial_k(u^k\Delta u^j)=\Delta\dot u^j
+\partial_i\big[(\operatorname{div} u)\partial_i u^j
-\partial_i u^k\partial_k u^j\big] -\partial_k(\partial_i u^k\partial_i u^j).
\end{align}

For the divergence term, one gets that
\begin{align*}
\operatorname{div} u_t=\operatorname{div}\dot u-\partial_i(u^k\partial_k u^i)
=\operatorname{div}\dot u-\partial_i u^k\partial_k u^i
-u^k\partial_k\operatorname{div} u,
\end{align*}
due to $u_t=\dot u-u\cdot\nabla u$. Hence
\begin{align}\label{c13}
\partial_j\operatorname{div} u_t=\partial_j\operatorname{div}\dot u
-\partial_j(\partial_i u^k\partial_k u^i)-\partial_j(u^k\partial_k\operatorname{div} u).
\end{align}
Moreover, one has
\begin{align*}
\partial_k(u^k\partial_j\operatorname{div} u)
&=\partial_k\big[\partial_j(u^k\operatorname{div} u)
-(\partial_j u^k)\operatorname{div} u\big]\notag\\
&=\partial_j\big[\partial_k(u^k\operatorname{div} u)\big]
-\partial_k\big[(\operatorname{div} u)\partial_j u^k\big]\notag\\
&=\partial_j\big[(\operatorname{div} u)^2+u^k\partial_k\operatorname{div} u\big]-\partial_k\big[(\operatorname{div} u)\partial_j u^k\big],
\end{align*}
which along with \eqref{c13} yields
\begin{align}\label{c15}
\partial_j\operatorname{div} u_t+\partial_k(u^k\partial_j\operatorname{div} u) =\partial_j\operatorname{div}\dot u+\partial_j\big[
(\operatorname{div} u)^2-\partial_i u^k\partial_k u^i\big] -\partial_k\big[(\operatorname{div} u)\partial_j u^k\big].
\end{align}

It follows from \eqref{eqpre} that
\begin{align}\label{c16}
-\partial_jP_t=\partial_j(u^k\partial_kP)+\gamma\partial_j(P\operatorname{div} u).
\end{align}
For the pressure transport term, a direct calculation gives
\begin{align*}
-\partial_k(u^k\partial_jP)&=-\partial_k\big[\partial_j(Pu^k)-P\partial_j u^k\big]\notag\\&=-\partial_j\big[\partial_k(Pu^k)\big]
+\partial_k(P\partial_j u^k)\notag\\&=-\partial_j(u^k\partial_kP)
-\partial_j(P\operatorname{div} u)+\partial_k(P\partial_j u^k),
\end{align*}
which combined with \eqref{c16} implies that
\begin{align}\label{c18}
-\partial_jP_t-\partial_k(u^k\partial_jP)=(\gamma-1)\partial_j(P\operatorname{div} u)+\partial_k(P\partial_j u^k).
\end{align}
Furthermore, the left-hand side of \eqref{c8} can be expanded as
\begin{align}\label{c9}
\partial_t(\rho\dot u^j)+\partial_k(u^k\rho\dot u^j)=\rho\partial_t\dot u^j
+\rho u^k\partial_k\dot u^j+\big[\rho_t+\partial_k(\rho u^k)\big]\dot u^j =\rho\left(\partial_t\dot u^j+u\cdot\nabla\dot u^j\right).
\end{align}
As a consequence, inserting \eqref{c12}, \eqref{c15}, \eqref{c18}, and \eqref{c9} into \eqref{c8}, we derive \eqref{c19}. \hfill$\square$

\end{document}